\documentclass{amsproc}
\usepackage{amsfonts}
\usepackage{cases}
\usepackage{mathrsfs}
\usepackage{bbm}
\usepackage{amssymb}
\usepackage{txfonts}
\usepackage{amscd}
\usepackage{amsfonts,latexsym,amsmath,amsthm,amsxtra,mathdots}
\usepackage[bookmarks,colorlinks]{hyperref}
\usepackage[all,cmtip]{xy}
\RequirePackage{amsmath} \RequirePackage{amssymb}
\usepackage{color}
\usepackage{colordvi}
\usepackage{multicol}
\usepackage{tikz}
\usetikzlibrary{matrix,arrows}

\newcommand{\alg}{{\text{alg}}}

\def\scrO{\mathscr{O}}

\def\an{\mathrm{an}}
\def\loc{\mathrm{loc}}

    \newcommand{\BA}{{\mathbb {A}}}

    \newcommand{\BG}{{\mathbb {G}}}

     \newcommand{\BN}{{\mathbb {N}}}
     \newcommand{\BP}{{\mathbb {P}}}
    \newcommand{\BQ}{{\mathbb {Q}}}

     \newcommand{\BZ}{{\mathbb {Z}}}

    \newcommand{\CI}{{\mathcal {I}}}

    \newcommand{\CO}{{\mathcal {O}}}

    \newcommand{\CU}{{\mathcal {U}}}

     \newcommand{\RB}{{\mathrm {B}}}
    \newcommand{\RC}{{\mathrm {C}}} 
     
    \newcommand{\RG}{{\mathrm {G}}} 
     
     \newcommand{\RL}{{\mathrm {L}}}
    \newcommand{\RM}{{\mathrm {M}}} 
     \newcommand{\RP}{{\mathrm {P}}}
     
     \newcommand{\RT}{{\mathrm {T}}}
    \newcommand{\RU}{{\mathrm {U}}}

    \newcommand{\Ad}{{\mathrm{Ad}}}

    \newcommand{\End}{{\mathrm{End}}} 
     \newcommand{\GL}{{\mathrm{GL}}}
     
    \newcommand{\Ind}{{\mathrm{Ind}}}

    \newcommand{\fo}{{\mathfrak{o}}}

     \newcommand{\rank}{{\mathrm{rank}}}

\def\bfC{\mathbf{C}}
\def\bfg{\mathbf{g}}
\def\bfG{\mathbf{G}}

\def\bfL{\mathbf{L}}

\def\bfP{\mathbf{P}}

\def\bfU{\mathbf{U}}
\def\bfu{\mathbf{u}}

\def\bfx{\mathbf{x}}

\def\bfOmega{\mathbf{\Omega}}

\def\fg{\mathfrak{g}}

\def\FH{\mathfrak{H}}
\def\FK{\mathfrak{K}}
\def\FU{\mathfrak{U}}

\def\scrO{\mathscr{O}}

\def\scrN{\mathscr{N}}

\def\pr{\mathrm{pr}}

\def\bfG{\mathbf{G}}

\def\-{^{-1}}

\newcommand{\delete}[1]{}

     \newcommand{\SL}{{\mathrm{SL}}}
     
     \newcommand{\Sym}{{\mathrm{Sym}}}

        \newcommand{\Sp}{{\mathrm{Sp}}}

    \newcommand{\ds}{\displaystyle}

    \newcommand{\ra}{\rightarrow}

    \theoremstyle{plain}

    \newtheorem{thm}{Theorem}[section] \newtheorem{cor}[thm]{Corollary}
    \newtheorem{lem}[thm]{Lemma}  \newtheorem{prop}[thm]{Proposition}
     \newtheorem{defn}[thm]{Definition}
     \newtheorem {example}[thm]{Example}
    
    \numberwithin{equation}{section}

\begin{document}

\title{Morita's duality for split reductive groups}
\author{Zhi Qi}

\begin{abstract}
In this paper, we extend the work in \emph{Morita's Theory for the
Symplectic Groups} \cite{QiYang} to split reductive groups. We
construct and study the holomorphic discrete series representation
and the principal series representation of a split reductive group
$\RG$ over a $p$-adic field $F$ as well as a duality between certain
sub-representations of these two representations.
\end{abstract}


\maketitle

{\footnotesize
\tableofcontents
}

\section*{Notations}

Let $p$ be a prime, $F$ a finite extension of $\BQ_p$, $\fo$ the
ring of integers of $F$, $\varpi$ a uniformizer of $\fo$, $|\ |$
the normalized absolute value, and $F^\alg$ an algebraic closure of
$F$. Let $K$ be an extension of $F$ with an absolute value extending
$|\ |$, and $\fo_K$ the valuation ring of $K$. We assume that $K$ is
complete with respect to $|\ |$, and moreover, $K$ is spherically complete whenever topological
properties of the $K$-vector spaces are under consideration.

\addtocontents{toc}{\protect\setcounter{tocdepth}{1}}

\section{Introduction}\label{sec:Intro}

In a series of   papers, \cite{MoritaI}, \cite{MoritaII} and
\cite{MoritaIII}, Morita and Murase innitiated the work on the
representation theory for $\SL(2, F)$ with coefficient field $K$, especially   
holomorphic discrete series representations and their duality
relations with   principal series representations. Principal
series representations, or more generally,   induced
representations, appeared in many literatures, notably in F\'eaux de
Lacroix's work \cite{Feaux} on   locally analytic representations.
On the other hand,   holomorphic discrete series representations were not
extensively studied. The holomorphic discrete series representation 
of $\SL(n+1, F)$ associated to a rational representation of $\GL(n,
F)$ were introduced by Schneider in \cite{Schneider1992}, in order
to understand the de Rham complex over Drinfel'd's space
as representation  of $\SL(n+1, F)$. In another
direction, the holomorphic discrete series representation  of
$\Sp(2n, F)$ associated to a $K$-rational representation of $\GL(n,
F)$ were constructed in our recent work \cite{QiYang}. Furthermore,
the algebraization and generalization of Morita's duality were
established in \cite{QiYang}.

The purpose of this paper is to generalize Morita's theory from
$\Sp(2n, F)$ to a split reductive group $\RG$. We are able to do such a
generalization owning to the entirely algebraic construction of
Morita's theory for $\Sp(2n, F)$. Therefore, we shall closely follow the
main ideas presented in \cite{QiYang}.

In the first paragraph, we recollect some notions on split reductive
groups and construct an $F$-regular function $f$ on $\RG$ that
characterizes the parabolic big cell associated to a parabolic
subgroup. In particular, $f$ corresponds to the determinant function on $\RM(n,
F)$ used in the definition of $p$-adic Siegel upper half-space in
\cite{QiYang} (see Example \ref{ex:symplectic} and \ref{ex:Siegel}). $f$ will appear extensively in the study of the rigid
symmetric space associated to $\RG$ and   holomorphic series
representations.

The principal series representation $(C^{\an}_\sigma({\FH} , V),
T_{\sigma})$ is another interpretation of the parabolic induction
from a locally analytic $K$-representation $(\sigma, V)$ of the Levi
subgroup (cf. \cite{MoritaII} and \cite{QiYang}). In the second
paragraph, applying the general results of F\'eaux de Lacroix on  
induced representations of $F$-Lie groups (\cite{Feaux}), one sees
that $(C^{\an}_\sigma({\FH} , V), T_{\sigma})$ is a locally analytic
representation of $\RG$ over a $K$-vector space of compact type.

The third paragraph is dedicated to the construction and study of
the rigid analytic symmetric space $\mathbf{\Omega}$, which is the
foundation of   holomorphic discrete series representations. Some
examples are the $p$-adic upper half-plane for $\SL(2, F)$ (cf.
\cite{MoritaI}), Drinfel'd's space for $\SL(n+1, F)$ (cf.
\cite{Schneider1992}) and the $p$-adic Siegel upper half-space for
$\Sp(2n, F)$ (cf. \cite{QiYang}). Such symmetric spaces have been studied
by van der Put and Voskuil in \cite{PutVoskuil}. We shall however
use another approach following \cite{Schneider-Stuhler} and
\cite{QiYang} to construct the admissible affinoid covering, which
enables us to obtain precise descriptions of   rigid analytic
functions on $\mathbf{\Omega}$. From this, we prove that the space
$\scrO_K(\bfOmega)$ of $K$-rigid analytic functions on $
\mathbf{\Omega}_K$ is a nuclear $K$-Fr\'echet space.

In the fourth paragraph, for a $K$-rational representation $(\sigma,
V)$ of the Levi subgroup, we construct the holomorphic discrete
series representation $(\mathscr{O}_\sigma(\Omega), \pi_{\sigma})$
defined over the nuclear $K$-Fr\'echet space of $V$-valued rigid
analytic functions on $\mathbf{\Omega}$. Moreover, we prove that its dual
representation is locally analytic.

Since the strong duality gives rise to a contra-variant equivalence between
the category of  $K$-vector spaces of compact type and the category of
nuclear $K$-Fr\'echet spaces (cf. \cite{Schneider-Teitelbaum2002}), 
it is natural to expect certain
duality relations between   sub-quotient spaces of   principal series
representations and those of   holomorphic discrete series
representations. For $\SL(2, F)$, a duality of this kind via   residues is analytically constructed by Morita
 (cf. \cite{MoritaII}). However, there does
not seem to be any direct way to generalize Morita's duality.
Nevertheless, 
Morita's duality may be algebraically interpreted
 in a weaker form and generalized  to any split
reductive group $\RG$. These are done in the last paragraph. For a
$K$-rational representation $(\sigma, V)$ of the Levi subgroup,  
 a closed sub-representation $B_{\sigma^*}(\FH,
V^*)$ of $C^{\an}_{\sigma^*}({\FH} , V^*)$ and a closed
sub-representation $\mathscr{N}_\sigma(\Omega)$ of
$\mathscr{O}_\sigma(\Omega)$ are algebraically constructed, along with a duality between them.
As discussed in \cite{QiYang}, our duality for $\SL(2, F)$ is exactly
  Morita's duality composed with Casselman's intertwining operator
(cf. \cite{MoritaII}).

\begin{small}
\emph{Acknowledgements.} The author is especially indebted to
Bingyong Xie for the communication of several important ideas.
Thanks are also due to Professor P. Schneider and Professor J.
Cogdell for comments and advices.
\end{small}

\section{A lemma on the split reductive groups}\label{sec:SRG}

\subsection{Split reductive groups}
We adopt the notations in \cite{Jantzen} Part II, Chapter 1.

 Let $\RG$ be a connected split reductive algebraic group over $F$, $\RT$ a
split maximal torus of $\RG$. We have the decomposition of Lie
algebra $\fg$ of $\RG$ (over $F$) in the form
\begin{equation}\label{eq:Liedecomposition}\fg = \fg_0
\oplus \bigoplus_{\alpha \in R} \fg_\alpha,\end{equation} where
$\fg_0$ is the Lie algebra of $\RT$ and $R$ is the root system of
$\RG$ with respect to $\RT$.

Each $\fg_\alpha$ is of rank $1$ over $F$, and we denote $\RU_\alpha
\simeq \BG_a(F)$ the root subgroup of $\RG$ corresponding to
$\alpha$.


Let $W \cong N_\RG (\RT) /\RT$ be the Weyl group of $R$. For
$w \in W$, we also denote $w$ a representing element in $N_\RG (\RT)$.

Choose a positive system $R^+$ and denote $S$ the corresponding set
of simple roots. Let $\RB^+ $ denote the corresponding Borel
subgroup and $\RB^- $ its opposite Borel subgroup, $\RU^\pm_\RG =
\RU(\pm R^+)$ the unipotent radical of $\RB^\pm$.

Throughout this article we fix a subset  $I$ of $S$, and denote $R_I
= R \cap \BZ I$, $W_I$ the Weyl group of $R_I$, $R^+_I = R^+ - R_I$,
$\RP^+$ the standard parabolic subgroup of $\RG$ corresponding to
$R^+_I$, $\RP^-$ its
opposite parabolic subgroup
, $\RU^\pm = \RU(\pm R^+_I)$ the unipotent radical of $\RP^\pm$,
$\RL $ the common Levi subgroup of $\RP^+ $ and $\RP^- $.

$\RL$ is a split reductive group with split maximal torus $\RT$,
root system $R_I$, positive system $R^+_\RL = R_I\cap R^+$ and Weyl
group $W_I.$ Let $\RU^\pm_\RL = \RU(\pm R^+_\RL).$

We recall (\cite{Jantzen} Part II.1.7) that for any closed and
unipotent subset $R'$ of $R$ (that is, $(\BN\alpha + \BN\beta) \cap R
\subset R'$ for any $\alpha, \beta \in R'$ and $R' \cap (-R') =
\emptyset $), for instance $\pm R^+$, $\pm R^+_I$ and $\pm R^+_\RL$,
the multiplication induces, for any ordering of $R'$, an isomorphism
of schemes over $F$
\begin{equation} \label{eq:U(R')=Ualpha}
\prod_{\alpha \in R' } \RU_\alpha \xrightarrow[\quad ]{\simeq } \RU
(R') .
\end{equation}

\subsection{The parabolic big cell}

We have the Bruhat decomposition of $\RG$ (\cite{Jantzen} Part
II.1.9)
$$\RG = \bigcup_{w \in W} \RB^- w \RB^+ = \bigcup_{w \in W} \RU^-_\RG w \RT \RU^+_\RG.$$
Let $\RC$ denote the parabolic big cell $$ \RC = \bigcup_{w \in W_I}
\RU^-_\RG w \RT \RU^+_\RG = \RU^-  \bigg( \bigcup_{w \in
W_I}\RU^-_\RL w \RT \RU^+_\RL \bigg) \RU^+ = \RU^- \RL \RU^+ = \RU^-
\RP^+ = \RP^- \RU^+.
$$
Then $\RG$ is the disjoint union of $\RC$ and $\RU^-_\RG w \RT
\RU^+_\RG$ for all $w \notin W_I$.

Let $ r = |R^+_I| = \dim \RU^+$. We consider the adjoint
representation of $\RG$ on $\bigwedge^r \fg$ over $F$. From the
decomposition (\ref{eq:Liedecomposition}) of $\fg$ we obtain a
direct sum decomposition of $\bigwedge^r \fg$. Choosing $X_\alpha$ a
nonzero element in $\fg_{\alpha }$ for each $\alpha \in R$ and a
basis of $\fg_0$ we obtain a basis of $\bigwedge^r \fg$ with respect
to this decomposition and containing $Y = \bigwedge_{\alpha \in
R^+_I} X_\alpha$. For $g \in \RG$ we define $f(g)$ to be the
coefficient of $Y$ in the expansion of $\Ad(g) Y$ in the chosen
basis of $\bigwedge^r \fg.$ Then $f$ is a regular function on $\RG$
over $F$.

There is a partial order on $\BZ S$: $ \gamma \prec \delta $ iff
$ \delta - \gamma $ is a sum of positive roots. If one considers the group action of the symmetric group 
$S_r$ on $( \BZ S)^r$ via coordinate permutation, the set of the unordered $r$-tuples $[\gamma_1, ..., \gamma_r]$
of elements in $\BZ S$ may be viewed  as the set of $S_r$-orbits
in $( \BZ S)^r$. We define $[\gamma_1, ..., \gamma_r] \prec [\delta_1, ..., \delta _r]$ iff
there exists $s \in S_r$ such that $\gamma_{j} \preceq \delta_{s(j)}$ for all $1 \leqslant  j \leqslant  r$
and $\gamma_{j} \prec \delta_{s(j)}$ for at least one $j$.

We adopt the convention that $\fg_\gamma = 0$ if $\gamma \in \BZ S$
is nonzero and not a root. Then $[\fg_\beta , \fg_\alpha] \subset
\fg_{\alpha + \beta}$, and therefore
$$\Ad(u_\beta) X_\alpha \in X_\alpha + \sum_{i \geqslant  1} \fg_{\alpha +
i \beta}, \hskip 10pt u_\beta \in \RU_\beta.$$

If $\alpha, \beta \in R_I^+$, then it is clear that either $\alpha +
i \beta \in R^+_I$ or $\fg_{\alpha + i \beta} = 0$. The same
statement holds for $\alpha \in R_I^+$ and $\beta \in R_I$, since
the $\beta$-string through $\alpha$ lies in $R^+_I$. Therefore
$U_\beta$ fixes $Y$ for any $\beta \in R_I^+ \cup R_I = R^+ \cup
(-R^+_I)$. (\ref{eq:U(R')=Ualpha}) implies that
\begin{equation}\label{eq:Y inv} Y \text{ is invariant under } \RU^+_\RG 
\text{ and }\RU^-_\RL. \hskip 188pt\end{equation}


If we let $\beta$ be negative roots, then (\ref{eq:U(R')=Ualpha})
also implies that for $\varv^- \in \RU^-_\RG $,
\begin{equation}\label{eq:Ad(U-G)} \Ad(\varv^-) X_\alpha \in X_\alpha +
\sum_{\gamma \prec \alpha} \fg_\gamma. \end{equation}

For $t \in T$, \begin{equation}\label{eq:Ad(t)} \Ad(t) Y =
\prod_{\alpha \in R^+_I} \alpha(t) Y.\end{equation}

For $w \in W$ there exists a constant $c_{w, \alpha} \in F^\times $
satisfying
\begin{equation}\label{eq:Ad(w)}\Ad(w) X_\alpha = c_{w, \alpha} X_{w
\alpha}.\end{equation}
 In view of (\ref{eq:U(R')=Ualpha}), we see that $w $ preserves $R^+_I$ iff $w
$ normalizes $\RU^+$. Since $N_\RG (\RU^+) = \RP^+$ and $w \in
\RP^+$ iff $w \in W_I$, 
\begin{equation}\label{eq:w preserves R+I} w \text{ preserves } R^+_I \text{ iff } w \in W_I. \hskip 210pt\end{equation}

Since $\RP^+ = \bigcup_{w \in W_I} \RU^-_\RL w \RT \RU^+_\RG$,
if we write $p^+ = u^- w t \varv^+$ ($\varv^+ \in \RU^+_\RG$, $\ t
\in \RT$, $ w \in W_I$, $u^- \in \RU_\RL^-$), then it follows from
(\ref{eq:Y inv}), (\ref{eq:Ad(t)}), (\ref{eq:Ad(w)}) and (\ref{eq:w
preserves R+I}) that
$$ \Ad (p^+ ) Y = \mathrm{sign}(w) \prod_{\alpha \in R^+_I} c_{w,
\alpha} \alpha(t)\cdot Y,
$$ where $ \mathrm{sign}(w)$ denotes the sign
of the permutation $w$ on $R^+_I.$ Moreover, it follows from
(\ref{eq:Ad(U-G)}) that
 $$\Ad (\varv^- w t \varv^+) Y \in ~ \mathrm{sign}(w) \prod_{\alpha \in R^+_I} c_{w,
\alpha} \alpha(t)\cdot Y + \sum_{[\gamma_j] \prec [\alpha]_{\alpha \in  R^+_I}} \bigwedge_{j = 1}^r \fg_{\gamma_j}.$$
So
$$ f(\varv^- w t \varv^+) = \mathrm{sign}(w) \prod_{\alpha \in R^+_I} c_{w,
\alpha} \alpha(t).
$$

Similarly, for $w \notin W_I$,
$$\Ad (\varv^- w t \varv^+) Y \in \sum_{[\gamma_j] \preceq [w \alpha]_{\alpha \in  R^+_I}} 
\bigwedge_{j = 1}^r \fg_{\gamma_j}.$$
It follows from the proof of (\ref{eq:Y inv}) that $\alpha + \beta \in R^+_I$ or
$\fg_{ \alpha + \beta } = 0$ if $\alpha \in R^+_I$ and $\beta \in R^+$, so $\alpha \in R^+_I$
and $ \delta \succeq \alpha $ imply $\delta \in  R^+_I$ or $\fg_\delta = 0$.
Therefore if $\delta \in \{0\} \cup R - R^+_I$ and $\gamma \preceq  \delta $ then $\gamma \notin R^+_I$,
and hence (\ref{eq:w preserves R+I}) implies that $Y$ does not appear in
the expression of $\Ad (\varv^- w t \varv^+) Y$, so $f(\varv^- w t \varv^+) = 0$.

We conclude with the following lemma.

\begin{lem} \label{lem:f(g)} Let the notations be as above, then

 $(1)$ For $p^+ \in \RP^+,$ $$\Ad(p^+) Y = f(p^+) Y,$$ and hence $f$
 is an $F$-rational character on $\RP^+$.

 $(2)$ For $w \in W_I,$
$$ f(\varv^- w t \varv^+) = \mathrm{sign}(w) \prod_{\alpha \in R^+_I} c_{w,
\alpha} \alpha(t), \hskip 10pt \varv^\pm \in \RU^\pm_\RG,\ t \in
\RT.
$$

 $(3)$ $f$ vanishes on $ \RU^-_\RG w \RT \RU^+_\RG$ for $w
\notin W_I.$\\
 In particular,

 $(4)$ $\RC$ is an open $F$-subscheme of $\RG$, and $ F[\RC] = F[\RG]_{f}.$

 $(5)$ $f$ is right invariant under $\RU^+_\RG$ and left invariant
under $\RU^-_\RG$.

\end{lem}

\begin{example} \label{ex:Borel} [cf. \cite{Jantzen} Part II 1.9 and \cite{Steinberg} \S 5 Theorem 7]
If $I = \emptyset$, then $\RL = \RT$ and $\RU^\pm = \RU^\pm_\RG$. Lemma
\ref{lem:f(g)} implies that $f(u^- t u^+) = \prod_{\alpha \in R^+ }
\alpha(t)$ and $ f(u^- w t u^+) = 0$ for all nontrivial $w$, $u^\pm
\in \RU^\pm_\RG$ and $ t \in \RT.$
\end{example}

\begin{example} \label{ex:SL(n+1)} Let
$\RG = \SL(n + 1, F)$.
Write $g = \begin{pmatrix}A &B
\\ C & d\end{pmatrix} $ with $A \in \RM(n, F), B \in \RM(n, 1 ; F), C \in \RM(1, n; F)$ and $d \in F$. Let
\begin{equation*}\begin{split}
\RU^+&  = \left\{\begin{pmatrix} I_n & 0 \\ C & 1
\end{pmatrix} \in \RG \right\},\\
\RL &  = \left\{\begin{pmatrix} A & 0 \\ 0 & d\end{pmatrix} \in \RG \right\}.
\end{split}\end{equation*}
 Some calculations show that $$f(g) = d^{n+1}.$$
\end{example}

\begin{example} \label{ex:GL(n)} More generally, we consider $\RG = \GL(n, F)$. Let $(n_1,..., n_s)$ be a partition
of $n$. Write $g = (g_{ij})_{1\leqslant  i, j\leqslant  s}$ with $g_{ij} \in
\RM(n_i, n_j; F)$. Let $\RU^+$ be the subgroup consisting of the
matrices $u = (u_{ij})$ such that $u_{ij} = 0$ for $j < i$ and
$u_{ii} = I_{n_i}$, and $\RL $ the subgroup consisting of the
matrices $l = (l_{ij})$ such that $l_{ij} = 0$ for $i \neq j$.
$\RL \cong \prod_{1\leqslant  i \leqslant  s} \GL(n_i, F)$.

For $l \in \RL$, \begin{equation*}\begin{split}
f(l) &= \prod_{1\leqslant  j < i
\leqslant  s} \det(l_{ii})^{n_j} \det(l_{jj})^{-n_i}\\
& = \prod_{1\leqslant  i \leqslant  s} \det(l_{ii})^{\sum_{j < i} n_j - \sum_{i
< k} n_k}.
\end{split}
\end{equation*}
The computation of the explicit formula for $f$ involves the process
of block lower (or upper) triangularization, and it turns out to be
complicated if $s \geqslant  3$. For $s = 2$,  we have
$$f(g) = \det(g_{22})^{n_1+n_2} \det(g)^{-n_2}.$$

The situation is the same if $\RG = \SL(n, F)$. For instance, if $s
= 2$, then
$$f(g) = \det(g_{22})^{n_1+n_2}.$$
\end{example}

\begin{example} \label{ex:symplectic} We consider $$\RG = \Sp(2n, F) = \left\{g \in \GL(2n, F)\ :\ {^tg} \begin{pmatrix}0 &I_n \\
-I_n & 0\end{pmatrix}g = \begin{pmatrix}0 &I_n \\
-I_n & 0\end{pmatrix} \right\}.$$ Write $g = \begin{pmatrix}A &B
\\ C & D\end{pmatrix} $ with $A, B, C, D \in \RM(n, F)$. Let
\begin{equation*}\begin{split}
\RU^+ &  = \left\{\begin{pmatrix}I_n & 0 \\ C & I_n \end{pmatrix}
\ :\ C \in \Sym(n, F) \right\},\\
\RL &   = \left\{\begin{pmatrix}{^t D\-} & 0 \\ 0 & D\end{pmatrix} \
:\ D \in \GL(n, F) \right\},
\end{split}\end{equation*}
where  $\Sym(n, F)$ is the group of symmetric matrices of order $n$ over $F$. Some calculations show that $$f(g) =
\det(D)^{n+1}.$$
\end{example}

\section{$\Ind_{\RP^-}^{\RG}\sigma$ and the principal series
$(C^{\an}_\sigma({\FH} , V), T_{\sigma})$}\label{sec:principal}

Induced representations, especially the parabolic inductions, are of
extreme importance in Lie theory. For $p$-adic Lie groups, they
were studied by F\'eaux de Lacroix in his work (\cite{Feaux}) on the
locally analytic representations.

We first recall the notion of locally analytic representations
over $K$ of an $F$-Lie group.

\begin{defn} [cf. \cite{Feaux} and \cite
{Schneider-Teitelbaum2002} \S 3] A \emph{locally analytic}
representation $(\sigma, V)$ of an $F$-Lie group $\RG$ on a
barreled locally convex Hausdorff $K$-vector space $V$ is a
continuous representation such that the orbit maps are $V$-valued
locally analytic functions on $\RG$. More precisely, for any $v \in
V$ there exists a BH-space $W$ of $V$ (that is, a Banach space $W$
together with a continuous injection $W \hookrightarrow V$) such
that $g \mapsto \sigma(g)v$ expands (in a neighborhood of the unit
element) to a power series with $W$-coefficients.
\end{defn}

Let $(\sigma, V)$ be a locally analytic representation of the Levi
subgroup $\RL$. $\sigma$ extends to a representation of $\RP^-$ defined by
$\sigma ( u l ) = \sigma(l)$ ($l \in \RL$, $u \in \RU^-$).

\begin{defn} Let $\Ind_{\RP^-}^\RG\sigma$ be the space of
$V$-valued locally analytic functions $\phi$ on $\RG$ satisfying $$
\phi(p^-g)=\sigma(p^-)\phi(g), \hskip 10pt \text{for all } g \in
\RG, p^- \in \RP^-. $$ On $\Ind_{\RP^-}^\RG\sigma$ we have a
$\RG$-action defined by right translation.
\end{defn}

Since the quotient space $ \RP^- \backslash \RG$ is compact, we
obtain from \cite{Feaux} 4.1.5 the following proposition.

\begin{prop}$\Ind_{\RP^-}^\RG\sigma$ is a locally analytic representation
of $\RG$. \end{prop}

Next, we introduce the principal series representation,   serving
as another description of $\Ind_{\RP^-}^\RG \sigma$.

Let $\FH$ and ${\overline{\FH}}$ denote the $\RG$-homogeneous spaces
$\RU^-\backslash \RG$ and $\RP^-\backslash \RG$ respectively, and
denote  $\hat{g} := \mathrm{pr}^\RG_{\FH}(g).$ Because $\RP^-\cong
\RU^-\rtimes \RL$, there is a left $\RL$-action on $\FH$, and
${\overline{\FH}} = \RL\backslash \FH$.

\begin{defn}Let $C^{\an}_\sigma({\FH} , V)$ be the space of $V$-valued locally
analytic functions $\varphi$ on ${\FH} $ satisfying
$$\varphi(l \hat{g})=\sigma(l) \varphi(\hat{g}), \hskip 10pt \text{ for all
} \hat{g}\in {\FH} \text{ and } l\in \RL.
$$
The \emph{principal series representation} $(C^{\an}_\sigma({\FH} ,
V), T_{\sigma})$ of $\RG$ is defined via
\begin{equation*}\label{eq:defTsigma}(T_{\sigma}(g) \varphi) (\widehat{g'})
:= \varphi(\widehat{g'} \cdot g).\end{equation*}

\end{defn}

\begin{lem} \label{lem:principalseries}~\\ $(1)$ The representation $\Ind_{\RP^-}^\RG\sigma$
is (naturally) isomorphic to $(C^\an_\sigma({\FH} , V), T_{\sigma})$.\\
$(2)$ $\Ind_{\RP^-}^\RG\sigma$ is isomorphic to
$C^{\an}({\overline{\FH}} , V)$.\\
$(3)$ Let $\iota$ be a locally analytic section of
$\mathrm{pr}^\FH_{{\overline{\FH}}}$, then $\iota$ induces an
isomorphism
 \begin{eqnarray*} 
 \iota^{\circ} :
C^{\an}_\sigma({\FH} , V) &\ra& C^{\an}({\overline{\FH}}
, V) \\
\nonumber \varphi &\mapsto& \varphi\circ \iota. 
\end{eqnarray*}
\end{lem}
\begin{proof}
(1) From a locally analytic section $\bar{\iota}$ of
$\mathrm{pr}^\RG_{{\FH}}$ we obtain an isomorphism (\cite{Feaux}
4.3.1) $$\bar{\iota}^\circ: \Ind_{\RU^-}^\RG\mathbf{1} \simeq
C^{\an} ({\FH} , V), \hskip 10pt \phi \mapsto \phi \circ
\bar{\iota}. $$ By restriction to the subspaces, $\bar{\iota}^\circ$ induces an
isomorphism, independent of $\bar{\iota}$, between $\Ind_{\RP^-}^\RG\sigma$ and
$C^{\an}_\sigma({\FH} , V)$.
$\RG$-equivariance is evident.

(2) A locally analytic section $\tilde{\iota}$ of
$\mathrm{pr}^\RG_{{\overline{\FH}}}$ induces an isomorphism
$\tilde{\iota}^\circ$ from $\Ind_{\RP^-}^\RG\sigma $ onto
$C^{\an}({\overline{\FH}} , V) $ (ibid.).

(3) Choose $\bar{\iota}$ and $\tilde{\iota}$  compatible with
$\iota$, that is, $\tilde{\iota} = \bar{\iota}\circ \iota$, then the
assertion follows from (1) and (2).
\end{proof}

Compactness of ${\overline{\FH}} $ implies that
$C^{\an}({\overline{\FH}} , V)$ is of compact type
(\cite{Schneider-Teitelbaum2002} Lemma 2.1). By
\cite{Schneider-Teitelbaum2002} Proposition 1.2, Theorem 1.3 and
\cite{Schneider} Proposition 16.10, we have the following corollary.

\begin{cor}\label{cor:compacttype} Suppose that $K$ is spherically complete. Let $ B$ be a closed subspace of
$C^{\an}_\sigma({\FH} , V)$, then both $B $ and $C^{\an}_\sigma({\FH} ,
V)/B$ are of compact type. In particular, they are reflexive,
bornological, and complete. Moreover, their strong duals $B^*_b$ and $(C^{\an}_\sigma({\FH},
V)/B)^*_b$ are nuclear Fr\'echet spaces.
\end{cor}

\section{Rigid analytic symmetric space $\bfOmega$}\label{sec:symmetricspace}

The rigid analytic symmetric space $\bfOmega$ associated to $\RG$
(with respect to a parabolic $\RP^+$) was constructed by van der Put
and Voskuil in \cite{PutVoskuil}. Some examples are the $p$-adic
upper half-plane, Drinfel'd's space and the $p$-adic Siegel upper
half-space, which are associated to $\SL(2, F)$, $\SL(n+1, F)$ and
$\Sp(2n, F)$ respectively (cf. \cite{MoritaI},
\cite{Schneider-Stuhler} and \cite{QiYang}).

\subsection{Definition of the symmetric space $\bfOmega$}

Let $\bfG$, $\bfP^\pm$, $\bfU^\pm$, $\bfL$ and $\bfC$ denote the
$F$-rigid analytifications of $\RG$, $\RP^\pm$, $\RU^\pm$, $\RL$ and
$\RC$ respectively. $f$ defines a rigid analytic function on $\bfG$.

Since $f$ is left invariant under $\RU^-$ (Lemma \ref{lem:f(g)} (5)), we
may define
$$f(\hat{g}, \bfu):= f(g \cdot \bfu)$$ for $\hat{g}\in \FH$ and
$\bfu \in \bfU^-$.

\begin{defn} \label{def:Omega} Let
\begin{equation*}\begin{split}
\bfOmega & :=\left\{ \bfu \in \bfU^- \ :\
g\cdot \bfu \in
\bfC, \text{ for all } g\in \RG \right\}\\
& \hskip 3.5pt = \left\{ \bfu \in \bfU^- \ :\ f(\hat{g}, \bfu) \neq 0, \text{ for
all } \hat{g}\in \FH \right\}.
\end{split}\end{equation*} We call $\bfOmega$ the
symmetric space associated to $\RG$ with respect to $\RP^+$.
\end{defn}

\begin{example}\label{ex:Drinfel'd} In the situation of Example
\ref{ex:SL(n+1)}, $\bfU^- \cong \mathbf{A}_{/F}^n$ and $(z_1, ...,
z_n) \in \bfOmega$ is given by the inequalities
$$c_1 z_1 + ... + c_n z_n + d \neq 0 \hskip 10pt \text{ for all nonzero } (c_1,...,c_n, d) \in F^{n+1}.$$
Therefore $\bfOmega$ is Drinfel'd's space.
\end{example}

\begin{example}\label{ex:Siegel} In the situation of Example
\ref{ex:symplectic}, $\bfU^- \cong \mathbf{Sym}(n)$ and $Z \in
\bfOmega$ is given by the inequalities
$$\det(C Z + D) \neq 0 \hskip 10pt \text{ for all } C\ {^tD}=D\ {^tC},\ \rank{ (C \; D)}=n.$$
Therefore $\bfOmega$ is the $p$-adic Siegel upper half-space.
\end{example}

We may also interpret $\bfOmega$ to be the complement of all the
$\RG$-translations of $(\bfG - \bfC)/ \bfP^+ = \bfG / \bfP^+ -
\bfU^- $ in $ \bfG / \bfP^+$. Therefore we have a left $\RG$-action
on $\bfOmega$ (induced from the left $\RG$-action on $\bfG /
\bfP^+$). We denote $g
* \bfu$ the action of $g\in \RG$ on $\bfu \in \bfOmega$. We have $g * \bfu = \pr^{\bfC}_{\bfU^-}(g \cdot
\bfu)$.

\subsection {Automorphy factor}

We define the \emph{automorphy factor }
\begin{eqnarray*} j : \RG \times \bfOmega &\ra&
\bfP^+ \\
(g, \bfu) & \mapsto & ( g * \bfu)^{-1}\cdot g\cdot \bfu.
\end{eqnarray*}
Then $j(g, \bfu) = \pr^{\bfC}_{\bfP^+}(g \cdot \bfu)$, and
straightforward computations show
\begin{equation}\label{eq:automorphyfactor} j(g_1  g_2, \bfu) = j(g_1,
g_2 * \bfu)j(g_2, \bfu).\end{equation}

For any $u \in \RU^- $, $j(u, \bfu) = 1_\bfG$, and hence
(\ref{eq:automorphyfactor}) implies $j(u\cdot g, \bfu) = j(g,
\bfu),$ so we may define $j(\hat{g}, \bfu) :=j(g, \bfu)$.

For $l \in \RL$, since $l * \bfu = l\cdot \bfu \cdot l\-$, we have
$j(l, \bfu) = l$, and by (\ref{eq:automorphyfactor})
\begin{equation}\label{eq:j(lg,bfu)}
j(l \cdot\hat{g}, \bfu) = l\cdot j(\hat{g}, \bfu).
\end{equation}

Since $f$ is left invariant under $\bfU^-$ (Lemma \ref{lem:f(g)}(5)),
it follows that
\begin{equation}\label{eq:f(j)}
f(j(\hat{g}, \bfu)) = f(\hat{g}, \bfu).
\end{equation}

From Lemma \ref{lem:f(g)} (1), (\ref{eq:automorphyfactor}) and
(\ref{eq:f(j)}), we see that
\begin{equation}\label{eq:f(g1g2&bfu)} f(\widehat{g_1 g_2}, \bfu) = f(\widehat{g_1}, g_2 *
\bfu) f(\widehat{g_2}, \bfu).
\end{equation}

Moreover, Lemma \ref{lem:f(g)} (1), (\ref{eq:j(lg,bfu)}) and
(\ref{eq:f(j)}) imply
\begin{equation}\label{eq:f(lg,bfu)} f(l \cdot \hat{g}, \bfu) = f(l) f(\hat{g}, \bfu).
\end{equation}

\subsection{The $F$-rigid analytic structure on $\bfOmega$}

van der Put and Voskuil defined an affinoid covering of $\bfOmega$
using Bruhat-Tits Buildings (\cite{PutVoskuil}). In this paper we
choose another approach following the construction of affinoid covering of Drinfel'd's space
in \cite{Schneider-Stuhler} and that of $p$-adic Siegel upper
half-space in \cite{QiYang}. We endow $\bfOmega $ with a structure
of $F$-rigid analytic variety and show that it is an admissible open
subset of $\bfU^-$ and, in particular, an open rigid analytic subspace
of $\bfU^-$ (and therefore $\bfG / \bfP^+$).

We realize $\RG$ as a Zariski closed subgroup of $\GL(n, F)$ such
that $\RT $ consists of diagonal matrices and $\RB^+ $ consists of
lower triangular matrices. Then $f(g)$ extends to an $F$-regular
function on $\GL(n, F)$ with respect to the coordinates $g_{i,j}$
($1\leqslant  i, j \leqslant  n$), and $\det(g)^r f(g)$ is a homogeneous
$F$-polynomial in $g_{i,j}$. We denote $N$ the degree of $\det(g)^r
f(g)$ and let $M$ be an integer such that all the coefficients have
absolute values bounded by $  |\varpi|^{NM}$.

\begin{lem}
$\bfOmega$ is nonempty.
\end{lem}
\begin{proof}
It suffices to prove that $\RG$-translations of $\bfG - \bfC$ do not
cover $\bfG$. For $g \in \RG$, $g \cdot (\bfG - \bfC)$ is the locus
of $f(g\- \cdot \bfg) = 0$. With the embedding of $\bfG$ into
$\mathbf{GL}(n, F)$ we view $f(g\- \cdot \bfg)$ as a rational
function in $\bfg_{ij}$ with $F[\RG]$-coefficients, and, for a given
$g \in \RG$, $f(g\- \cdot \bfg)$ is a nonzero $F$-rational function
in $\bfg_{ij}$. It is not hard to see that there are choices of
$\bfg_{ij} \in F^\alg $ with appropriate absolute
values so that the non-vanishing monomials in $f(g\- \cdot \bfg)$
are of distinct absolute values in $|(F^\alg)^\times| =
|\varpi|^{\BQ} $ modulo $|F^\times| = |\varpi|^{\BZ}$. Therefore there
exists $\bfg \in \bfG$ such that $f(g\- \cdot \bfg)$ is nonzero for
all $g \in \RG$, and consequently $\bfg$ lies in the complement of all the $g \cdot
(\bfG - \bfC)$.
\end{proof}

Let $\RG_\fo$ and $\RL_\fo$ denote the intersections of $\RG$ and
$\RL$ with $\GL(n, \fo)$ respectively, and denote $\FH_\fo =
\pr^\RG_\FH (\RG_\fo)$.

We recall Iwasawa's decomposition
$$\RG = \RB^- \RG_\fo.$$
Then $\RG = \RP^-\cdot \RG_{\fo}$ and $\FH=\RL \cdot \FH_\fo$, so
(\ref{eq:f(lg,bfu)}) and Lemma \ref{lem:f(g)} imply
$$ \bfOmega = \left\{ \bfu \in \bfU^- \ :\ f(\hat{g}, \bfu) \neq 0, \text{ for any } \hat{g}\in
\FH_\fo \right\}.
$$

For $\bfu\in \bfU^-$ an upper triangular matrix with diagonal
entries $1$, let
$$|\bfu| := \max_{1\leqslant  i \leqslant  j \leqslant  n}
 |\bfu_{ij}| = \max_{1\leqslant  i < j \leqslant  n} \left\{ 1,
 |\bfu_{ij}| \right\} .$$  
For any nonnegative integer $m$ and $\hat{g}
\in\FH_\fo $, we define
$$ \mathbf{B} (m; \hat{g}) := \left\{ \bfu\in \bfU^- \ :\ | f(\hat{g}, \bfu) |
< |\bfu|^{N }\ |\varpi|^{N (M + m)}\right\}.$$

\begin{lem} \label{lem:B(m;g)}
If $m$ is a nonnegative integer and $g_1, g_2\in\RG_\fo $ such that
$g_1 \equiv l\cdot g_2 \mod \varpi^{N m+1}$ for some $l\in \RL_\fo$,
then
$$\mathbf{B} (m;\widehat{g_1}) =\mathbf{B} (m; \widehat{g_2}).$$
\end{lem}

\begin{proof} Since $f|_{\RL_\fo}$ is a continuous $F$-character (Lemma
\ref{lem:f(g)} (1)) and $\RL_\fo$ is compact, the image of $\RL_\fo$
under $f$ is contained in $\fo^\times$. Therefore
(\ref{eq:f(lg,bfu)}) implies $|f(\widehat{l g_2}, \bfu)| =
|f(\widehat{g_2}, \bfu)|$, and hence $\mathbf{B} (m; \widehat{g_2}) =
\mathbf{B} (m; \widehat{l g_2}) $. So we may assume $ g_1 \equiv g_2
\mod \varpi^{N m+1}$.

We choose $\lambda \in (F^\alg)^\times$ such that
$|\lambda|=|\bfu|$. Since $|\lambda^{-1} \bfu_{ij}|\leqslant  1$,
$$ g_1\cdot {\lambda} \- \bfu \equiv
g_2\cdot {\lambda} \-  \bfu \mod \varpi^{N m+1},$$ and the matrices
on both sides have entries with absolute values $\leqslant  1$. Applying
$\det^r \cdot f$, we obtain
$${\lambda^{-N}} \det(g_1 )^r f(g_1\cdot \bfu) \equiv
{\lambda^{-N}} \det(g_2 )^r f(g_2\cdot \bfu) \mod \varpi^{NM + N
m+1},$$ and consequently
$$ |\bfu|^{-N}  \left|f(\widehat{g_1}, \bfu)\right| < |\varpi|^{N (M + m)}
\Leftrightarrow |\bfu|^{-N}  \left|f(\widehat{g_2}, \bfu)\right| <
|\varpi|^{N (M + m)}.
$$
Therefore $\mathbf{B} (m;\widehat{g_1}) =\mathbf{B} (m;
\widehat{g_2}) $.
\end{proof}

Let
\begin{equation*}\begin{split}
\mathbf{\Omega}(m; \hat{g}) &:=\bfU^- -  \mathbf{B} (m; \hat{g})\\
&\hskip 3pt = \left\{ \bfu \in \bfU^- \ :\ \left|f(\hat{g}, \bfu)\right| \geqslant 
|\bfu_{ij}|^{N } |\varpi|^{N (M + m)}, 1 \leqslant  i \leqslant  j \leqslant  n \right\}, \\
\mathbf{\Omega} (m) &:= \bigcap_{\hat{g}\in\FH_\fo} \mathbf{\Omega}
(m; \hat{g}).
\end{split}
\end{equation*}

For a given $\bfu \in \bfOmega$, $|f(\hat{g}, \bfu)|$ has a
positive lower bound on $\FH_\fo$. Therefore
$$\bfOmega = \ds \bigcup_{m =0}^\infty \mathbf{\Omega}(m).$$

Let $\FH^{(m)}$ be any finite subset of  $\FH_\fo$ including a set
of representatives in $\FH_\fo$ for $\pr^\RG_\FH \left(\RL_\fo \big \backslash \RG_\fo \big/ \RG_\fo(Nm+1) \right),$ where
$\RG_\fo (Nm + 1)$ denotes the congruence subgroup $\big(I_n + \varpi^{Nm+1} \RM(n, \fo) \big) \cap \RG $. Then Lemma \ref{lem:B(m;g)} implies that
$$\mathbf{\Omega}(m) = \bigcap_{\hat{g}\in\FH^{(m)}} \mathbf{\Omega} (m; \hat{g}).$$
Moreover, we may assume that $\FH^{(m)}$ contains $\hat{I_n}$.

$$\mathbf{\Omega}(m; \hat{I_n}) = \left\{\bfu\in \bfU^- \
:\ \left|\varpi ^{M + m} \bfu_{ij} \right| \leqslant  1 \right\}$$ is an admissible
open affinoid subset of $\bfU^-$. $\mathbf{\Omega}(m)$ is the
intersection of finitely many rational sub-domains of
$\mathbf{\Omega}(m; \hat{I_n} ) $:
$$\left\{\bfu\in \mathbf{\Omega}(m; \hat{I_n})\ :\
\left|\frac{\varpi^{N(M + m)} \bfu_{ij} ^{N}}{f(\hat{g},
\bfu)}\right|\leqslant  1, 1 \leqslant  i \leqslant  j \leqslant  n \right\},
$$ with $\hat{g}$ ranging on $\FH^{(m)} - \{\hat{I_n}\}$. Therefore $\mathbf{\Omega}(m) $ is an affinoid
variety.

We conclude that $\{\mathbf{\Omega}(m)\}_{m=0}^{\infty}$ constitutes an
admissible affinoid covering of $\mathbf{\Omega}$ so that
$\mathbf{\Omega}$ admits a rigid analytic variety structure (see
\cite{Bosch1984} 9.3). According to \cite{Bosch1984} 9.1.2 Lemma 3
(compare \cite{Bosch1984} 9.1.4 Proposition 2), the following
proposition implies that $\mathbf{\Omega}$ is an admissible open
subset of $\bfU^-$.

\begin{prop}\label{prop:Omegaadmissible}
Any morphism from an affinoid variety to $\bfU^-$ with image in
$\mathbf{\Omega}$ factors through some $\mathbf{\Omega}(m)$.
\end{prop}

\begin{proof}
 The argument is similar to the third proof of
\cite{Schneider-Stuhler} \S 1 Proposition 1.

Let $\mathbf{X}$ be an affinoid variety, $\Phi:
\mathbf{X}\rightarrow \bfU^- $ a morphism from $\mathbf{X}$ to
$\bfU^-$ with image in $\mathbf{\Omega}$. For any $\hat{g}
\in\FH_\fo$,
$$ \bfx \mapsto \frac{ \Phi(\bfx)_{ij} ^N}{f(\hat{g}, \Phi(\bfx))}, \hskip
10pt 1 \leqslant  i \leqslant  j \leqslant  n,
$$
are $F$-rigid analytic functions on $\mathbf{X}$. By the maximum
modulus principle (\cite{Bosch1984} \S 6.2 Proposition 4 (i)), there
exists a positive integer $m_{\hat{g}}$ such that
$$\max_{1\leqslant  i \leqslant  j \leqslant  n}
\max_{\bfx\in \mathbf{X}} \left|\frac{ \Phi(\bfx)_{ij} ^N}
 {f(\hat{g}, \Phi(\bfx))}\right|  \leqslant  |\varpi|^{- N (M + m_{\hat{g}})}.$$
In other words, $\Phi(\mathbf{X})\subset \mathbf{\Omega}(m_{\hat{g}};
\hat{g})$. In view of Lemma \ref{lem:B(m;g)}, $m_{\hat{g}}$ can be
chosen locally constant. Therefore the compactness of $\FH_\fo$
implies that there exists a positive integer $m$ such that
$\Phi(\mathbf{X}) \subset \mathbf{\Omega}(m)$.
\end{proof}

Finally, we prove that the morphisms of $g$-translations from
$\bfOmega(m)$ into $\bfOmega$ indeed factor through the same
$\mathbf{\Omega}(m')$ for all $g \in \RG_\fo$.

\begin{lem} \label{lem:gOmegam} For any nonnegative integer $m$, there exists a 
nonnegative integer $m'$ such that for all $g\in \RG_{\fo}$,
$$g * \mathbf{\Omega}(m) \subset \mathbf{\Omega}(m').$$
\end{lem}
\begin{proof}
Let $\bfu\in \bfOmega(m)$. Then
\begin{equation}
\label{eq:gOmega1} 1 \leqslant  |\bfu| \leqslant  |\varpi|^{ - M - m},
\end{equation}
and
\begin{equation}
\label{eq:gOmega2}\frac {|\bfu|^N} {\left|f(\hat{g}, \bfu)\right|}
\leqslant  |\varpi|^{- N(M + m)}  \hskip 10pt \text{for any } g
\in\RG_\fo.
\end{equation}

$g * \bfu = \pr^{\bfC}_{\bfU^-}(g \cdot \bfu)$, and since
$\pr^{\bfC}_{\bfU^-}$ is $F$-regular on $\bfC$, Lemma \ref{lem:f(g)}
(4) implies that there exist positive integers $s$ and $t$ such
that all the entries of $$\det(g \cdot \bfu)^t f(g \cdot \bfu)^s
\cdot g* \bfu = \det(g )^t f(g \cdot \bfu)^s \cdot g* \bfu$$ are
$F$-polynomials with variables the entries of $g \cdot \bfu.$ Let
$D$ be the highest degree and $L$ an integer such that the absolute
values of all the coefficients are bounded by $|\varpi|^{L}$. Since $g
\in\RG_\fo$, the
entries of $g \cdot \bfu$ have absolute values $\leqslant  |\bfu|$ and
$|\det(g )| = 1$, then
\begin{equation}
\label{eq:gOmega3}|g * \bfu| \leqslant  \frac {|\varpi|^{ L} |\bfu|^D } {
\left| f(\hat{g}, \bfu)\right|^s }.
\end{equation}


It follows from (\ref{eq:f(g1g2&bfu)}), (\ref{eq:gOmega3}),
(\ref{eq:gOmega1}) and  (\ref{eq:gOmega2}) that
for any $g_1 \in \RG_\fo$,
\begin{align*} 
\frac {|g * \bfu|^N}{|f(\widehat{g_1}, g * \bfu)|}  &\leqslant  \frac {|\varpi|^{N L} |\bfu|^{N D} }
{|f(\widehat{g_1 g}, \bfu)| | f(\hat{g}, \bfu)|^{Ns - 1} }\\
& =  |\varpi|^{N L} |\bfu|^{N D - N^2 s} \frac { |\bfu|^{N } }
{|f(\widehat{g_1 g}, \bfu)| }  \frac { |\bfu|^{N^2 s - N} } {
|f(\hat{g}, \bfu)|^{Ns - 1} },\\
& \leqslant   |\varpi|^{N(L - \max\{D, Ns\} (M + m))}.
\end{align*}
 Therefore $g * \mathbf{\Omega}(m) \subset \mathbf{\Omega}(m') $ for any $m' \geqslant  - M - L + \max\{D, Ns\} (M +m)$.
\end{proof}

\subsection{Rigid analytic functions on $\bfOmega$}

Let $\scrO(\mathbf{\Omega}(m))$ denote the space of $F$-rigid
analytic functions on $\mathbf{\Omega}(m)$. Then
$\scrO(\mathbf{\Omega}(m))$ is an $F$-affinoid algebra with the
supremum norm.

Let $\scrO(\mathbf{\Omega})$ be the $F$-algebra of $F$-rigid
analytic functions on $\mathbf{\Omega}$, that is, the projective
limit of $\scrO (\mathbf{\Omega}(m))$,
$$\scrO(\mathbf{\Omega}) :=\underset{m}{\varprojlim}\scrO (\mathbf{\Omega}(m)).$$ 
$\scrO(\mathbf{\Omega})$ is endowed with the projective limit topology.

From the construction of $\bfOmega (m)$ we see that the $F$-affinoid
algebra $\scrO (\mathbf{\Omega}(m))$ is equal to
 \begin{equation} \label{eq:generatorsofOmega(m)} F \left\langle \varpi^{M + m}
\bfu_{ij} , \frac{\varpi^{N (M + m)} \bfu_{ij}^N }{f(\hat{g},
\bfu)}\ :\  1\leqslant  i \leqslant  j \leqslant  n,\ \hat{g} \in \FH ^{(m)} -
\{\hat{I_n}\}\right\rangle.\end{equation}

Therefore $\psi \in \scrO(\mathbf{\Omega}(m))$ has an expansion in
the following form that converges with respect to the supremum norm $||\
||_{\scrO(\mathbf{\Omega}(m))}$:
\begin{equation}\label{eq:psibfuexpansion} \psi (\bfu) = \sum _{(\ell_{\hat{g}})\in
(\BN_0)^{\FH^{(m)}}} P_{(\ell_{\hat{g}})}(\bfu)\prod_{{\hat{g}}\in
\FH^{(m)}} f({\hat{g}, \bfu }) ^{-\ell_{\hat{g}}},\end{equation}
where $P_{(\ell_{\hat{g}})}(\bfu)$ are polynomials in the
coordinates of $\bfu$ with coefficients in $F$.

In view of (\ref{eq:f(lg,bfu)}), the assumption that
$\FH^{(m)}$ is contained in $\FH_\fo$ is quite artificial, and it is
more convenient and natural to choose $\FH^{(m)}$ to be an arbitrary
finite subset of $\FH$ whenever we consider the
expansion of $ \psi \in \scrO(\mathbf{\Omega}(m))$.

$\psi \in  \scrO(\mathbf{\Omega} ) $ may be considered as an
$F^\alg$-valued function on $\mathbf{\Omega}$ such that, restricting on each
$\bfOmega (m)$, $\psi$ has an expansion (\ref{eq:psibfuexpansion})
that converges with respect to $||\
||_{\scrO(\mathbf{\Omega}(m))}$. In particular, $f({\hat{g}, \bfu
})^{-1} \in \scrO(\mathbf{\Omega} )$ for any ${\hat{g}} \in \FH$.

Since all the generators of
$\mathscr{O}(\bfOmega(m))$ in  (\ref{eq:generatorsofOmega(m)}) are $F$-rigid analytic functions on
$\bfOmega(m')$ for any $m' \geqslant  m$ and therefore on $\bfOmega$, we
obtain the following proposition.

\begin{prop} \label{prop:Stein} \ \\
$(1)$ $\bfOmega$ is a Stein space, that is, the image
of $\mathscr{O}(\bfOmega(m + 1))$ under the transition homomorphism
in $\mathscr{O}(\bfOmega(m))$ is dense for any nonnegative integer
$m$.
\\
$(2)$ The image of $\mathscr{O}(\bfOmega)$ under the transition
homomorphism in $\mathscr{O}(\bfOmega(m))$ is dense.
\end{prop}

Let $\scrO_K(\mathbf{\Omega}(m))$ and $\scrO_K(\mathbf{\Omega})$
denote $\ds \scrO(\mathbf{\Omega}(m)) {\hat{\otimes}}_F  K $ and
$\scrO(\mathbf{\Omega}) {\hat{\otimes}}_F K$ respectively. If we let
$\mathbf{\Omega}_K(m)$ and $\mathbf{\Omega}_K $ denote the extensions
of the ground field $K/F$ of $\mathbf{\Omega} (m)$ and
$\mathbf{\Omega} $ respectively (see \cite{Bosch1984} \S 9.3.6), then
$\scrO_K(\mathbf{\Omega}(m))$ and $\scrO_K(\mathbf{\Omega})$ are the spaces of
$K$-rigid analytic functions on $\mathbf{\Omega}_K(m)$ and
$\mathbf{\Omega}_K$ respectively.

\begin{prop} \label{prop:OOmega} Let $K$ be spherically complete.
$\scrO_K(\mathbf{\Omega})$ is a nuclear $K$-Fr\'echet space.
\end{prop}
\begin{proof}
By \cite {Schneider} Proposition 19.9, it
suffices to prove that all the $\scrO_K(\mathbf{\Omega}(m))$ constitute a compact
projective system.

Consider the $\scrO_K(\mathbf{\Omega}(m-1))$-norms of the generators of
$\scrO_K(\mathbf{\Omega}(m))$ (see \ref{eq:generatorsofOmega(m)}), then
$$\sup_{\bfu\in \mathbf{\Omega}(m-1)} \max_{ \hat{g}\in
\FH ^{(m)} - \{\hat{I_n}\}} \max_{1\leqslant  i \leqslant  j \leqslant  n} \left\{
\left|\varpi^{M + m} \bfu_{ij} \right|, \left| \frac{\varpi^{N (M + m)}
\bfu_{ij} }{f(\hat{g}, \bfu)}\right| \right\}\leqslant  |\varpi|.$$
 \cite {Schneider-Teitelbaum2002BV} Lemma
1.5 implies that the transition homomorphism from $\scrO_K
(\mathbf{\Omega}(m))$ to $ \scrO_K (\mathbf{\Omega}(m-1))$ is
compact.
\end{proof}

\cite{Schneider-Teitelbaum2002}
Theorem 1.3 and Proposition 1.2 imply the following corollary.

\begin{cor}
Suppose $K$ is spherically complete. Let $\scrN$ be a closed subspace of $\scrO_K(\mathbf{\Omega})$, then
$\scrN$ and $\scrO_K(\mathbf{\Omega})/\scrN$ are nuclear Fr\'echet spaces, and their strong duals $\scrN^*_b$ and $( \scrO_K(\mathbf{\Omega})/\scrN )^*_b$ are of compact
type.
\end{cor}

\section{Holomorphic discrete series $(\mathscr{O}_\sigma(\Omega), \pi_{\sigma})$}
\label{sec:discrete}

Let $\Omega$ (resp.  $\Omega(m)$) denote $\mathbf{\Omega}_K (K)$
(resp. $\mathbf{\Omega}_K(m)(K)$). Restricting to $\Omega$ (resp.
$\Omega(m)$), we view $K$-rigid analytic functions in
$\mathscr{O}_K(\mathbf{\Omega} )$ (resp.
$\mathscr{O}_K(\mathbf{\Omega}(m))$) as $K$-valued functions on
$\Omega$ (resp. $\Omega(m)$), and abbreviate
$\mathscr{O}_K(\mathbf{\Omega}(m))$ (resp.
$\mathscr{O}_K(\mathbf{\Omega})$) to $\mathscr{O}(\Omega(m))$ (resp.
$\mathscr{O}(\Omega)$).

Let $(V,\sigma)$ be a $d$-dimensional $K$-rational representation of
$\RL$. $\sigma$ extends to a representation of $\RP^+$.

 Let $\mathscr{O}_\sigma(\Omega) :=\mathscr{O}(\Omega)\underset{
K}{\otimes} V$ and $\mathscr{O}_\sigma(\Omega(m))
:=\mathscr{O}(\Omega(m)) \underset{ K}{\otimes} V$. 

For any $ g \in \RG$ and $ \psi\in \mathscr{O}_\sigma(\Omega),$ let
$\pi_{\sigma}(g) \psi$ be the $V$-valued function on $\Omega$ as
follows
\begin{equation*}\label{eq:discrete}(\pi_{\sigma}(g) \psi) (\bfu) :=\sigma(j(g^{-1},
\bfu))^{-1}\psi(g^{-1}*\bfu).
\end{equation*}

\begin{lem}\label{lem:pi(psi)inO}
$\pi_{\sigma}(g) \psi \in \mathscr{O}_\sigma (\Omega).$
\end{lem}
\begin{proof}
 Since $\sigma$ is $K$-rational, each coordinate of $\sigma(j(g\-,
\bfu))^{-1} $ is a product of a $K$-polynomial in the coordinates of $j(g\-, \bfu)
$ and a power of $\det(j(g\-, \bfu))^{-1} =
\det(g) $. Note that $j(g\-, \bfu) = \pr^{\bfC}_{\bfP^+}(g\- \cdot \bfu)$,
and since $ \pr^{\bfC}_{\bfP^+}$ is $F$-regular on $\bfC$, each
coordinate of $j(g\-, \bfu) $ is a product of an
$F$-polynomial in the coordinates of $g\- \cdot \bfu$ and powers of $\det(g\-
\cdot \bfu)\- = \det(g)$ and $f(g\- \cdot \bfu)\-$. Therefore
each coordinate of $\sigma(j(g\-, \bfu))^{-1} $ has a finite
expansion of the form (\ref{eq:psibfuexpansion}), and hence belongs to
$\mathscr{O}(\Omega)$.

 Similarly, the coordinates of $\psi(g^{-1}*\bfu)$ also have
 expansions of the form (\ref{eq:psibfuexpansion}). By Proposition
 \ref{prop:Omegaadmissible}, for any $m$, $g\-$-translation maps $\bfOmega(m)$
 into some $\bfOmega(m')$, and hence the norm of each coordinate of $\psi(g^{-1}*\bfu)$ on $\bfOmega(m)$ is
 bounded by the norm of the corresponding coordinate of $\psi$ on $\bfOmega(m')$. Therefore $\psi(g^{-1}*\bfu) \in  \mathscr{O}_\sigma (\Omega).$

We conclude that $\pi_{\sigma}(g) \psi \in
\mathscr{O}_\sigma (\Omega).$
\end{proof}

It follows from the automorphy relation (\ref{eq:automorphyfactor})
that $\pi_\sigma$ is an action of $\RG$ on
$\mathscr{O}_\sigma(\Omega)$.

\begin{defn}
We call $(\mathscr{O}_\sigma(\Omega), \pi_{\sigma})$ the
\emph{holomorphic (rigid analytic) discrete series representation}
of $\RG$.
\end{defn}

\begin{lem}\label{lem:norms-m-m'} Let $m$ and $m'$ be as in Lemma
\ref{lem:gOmegam}. Then there exists a constant $c$ depending on
$\sigma$ and $m$ such that $$||\pi_{\sigma}(g)
\psi||_{\mathscr{O}_\sigma(\Omega(m))} \leqslant  c
||\psi||_{\mathscr{O}_\sigma(\Omega(m'))},$$ 
for all $g \in \RG_\fo.$
\end{lem}

\begin{proof}
The proof is similar to the arguments in Lemma \ref{lem:pi(psi)inO},
but instead of Proposition \ref{prop:Omegaadmissible} we apply Lemma
\ref{lem:gOmegam}.

Using the expressions for the coordinates of $\sigma(j(g\-, \bfu))^{-1} $ 
in the first paragraph of the proof of Lemma \ref{lem:pi(psi)inO}, we see that
their $\mathscr{O} (\Omega(m))$-norms are uniformly bounded on
$\RG_{\fo}$, so there is a constant $c > 0$ such that
$$\max_{g\in \RG_{\fo}} \max_{ \bfu\in
\mathbf{\Omega}(m)}||\sigma(j(g^{-1}, \bfu))^{-1}||_{\End(V)}\leqslant  c.
$$
Consequently,
\begin{align*} 
& \hskip 12.5pt \max_{g\in \RG_{\fo}} ||\pi_{\sigma}(g) \psi||_{\mathscr{O}_\sigma(\Omega(m))}\\
&= \max_{g\in \RG_{\fo}} \max_{\bfu\in \mathbf{\Omega}(m) } ||(\pi_{\sigma}(g) \psi)(\bfu)||_V \\
&\leqslant  \max_{g\in \RG_{\fo}}\max_{\bfu\in\mathbf{\Omega}(m)
} ||\sigma(j(g^{-1}, \bfu))^{-1}||_{\End(V)} \cdot
\max_{g\in \RG_{\fo}} \max_{\bfu\in \mathbf{\Omega}(m) } ||\psi(g^{-1} * \bfu)||_V \\
&\leqslant  c \max_{\bfu \in \mathbf{\Omega}(m') } ||\psi(\bfu)||_V \\
&= c ||\psi||_{\mathscr{O}_\sigma(\Omega(m'))}.
\end{align*}
\end{proof}

It follows from Lemma \ref{lem:norms-m-m'} that, for each $m$, the
map
\begin{eqnarray*} \RG_{\fo} \times \mathscr{O}_\sigma(\Omega )  & \rightarrow &
\mathscr{O}_\sigma(\Omega(m))  \\
\nonumber (g,\psi) & \mapsto & (\pi_{\sigma}(g) \psi)|_{\Omega(m)}.
\end{eqnarray*}
is continuous. Since $\mathscr{O}_\sigma(\Omega)$ is the projective
limit of $\mathscr{O}_\sigma(\Omega(m))$, we obtain the following corollary.

\begin{cor}\label{cor:picontinuous}
$(\mathscr{O}_\sigma(\Omega), \pi_{\sigma})$ is a continuous $\RG$-representation.
\end{cor}

Moreover, we shall prove that the dual representation of $\pi_\sigma$ is
locally analytic. For this, we recall that a coordinate chart at
$1_\RG$ is obtained from the decomposition of the Bruhat big cell
(see (\ref{eq:U(R')=Ualpha}))
 \begin{equation}\label{eq:coordinate} \RU^- \RU_\RL^- \RT
\RU_\RL^+ \RU^+ \simeq \BA_F^{|R|} \times
\BG_{m}(F)^{\dim \fg_0 }.\end{equation}

\begin{lem}\label{lem:pi_sigma analytic} Let $m$ and $m'$ be as in Lemma
\ref{lem:gOmegam}. Let $\RB  $ be any parameterized (as in
(\ref{eq:coordinate})) open neighborhood of $1_\RG$ contained in
$\RG_\fo$. For any $\psi\in \mathscr{O}_\sigma(\Omega(m'))$, the
orbit map
\begin{eqnarray*}  \RB & \rightarrow & \mathscr{O}_\sigma(\Omega(m))  \\
g & \mapsto & (\pi_{\sigma}(g) \psi)|_{\Omega(m)}
\end{eqnarray*} is an $\mathscr{O}_\sigma(\Omega(m))$-valued analytic function, namely, it can be expanded as a convergent power series with variables the
coordinate parameters of $\RB $ and coefficients in the Banach space
$\mathscr{O}_\sigma(\Omega(m))$.
\end{lem}

\begin{proof}
Once we have obtained a formal expansion of $\pi_\sigma (g)\psi$ into a
power series with variables the coordinate parameters of $\RB $ and
coefficients in $\mathscr{O}_\sigma(\Omega(m))$, Lemma
\ref{lem:norms-m-m'} would imply that the expansion is indeed
convergent. In view of (\ref{eq:coordinate}), it suffices to
consider $\pi_\sigma (g) \psi(\bfu) $ for $g $ in $ \RU^-,
\RU^-_\RL$ and $\RT$ (note that $\RU^+$ and $\RU^+_\RL$ are the
conjugations of $\RU^-$ and $\RU^-_\RL$ by the long Weyl element).

Let $u \in \RU ^-$, then
$$\pi_\sigma (u) \psi(\bfu) = \psi(u\- \cdot \bfu).$$

 Let $\mathscr{O}(\Omega(m))[[u]]$ denote the ring of formal power series $\varphi(u)$
 in the coordinates $u_{\alpha}$ ($\alpha \in R^-_I$) with coefficients in $\mathscr{O} (\Omega(m))$, where $\varphi(u)$ is expressed as
$$\varphi(u)=\sum_{\underline{r}\in \BN_0^{ R^-_I }} a_{{\underline{r}}}\cdot
 {\underline{u}}^{{\underline{r}}} , \hskip 10pt
a_{{\underline{r}}} \in \mathscr{O}(\Omega(m)),\
{\underline{u}}^{\underline{r}}:=\prod_{\alpha \in R^-_I}
u_{\alpha}^{r_{\alpha}}.$$
If the constant term
$a_{\underline{0}}$ is a unit in $\mathscr{O}(\Omega(m))$, then
$\varphi(u) $ is invertible in $\mathscr{O}(\Omega(m))[[u]] $. Note
that, for $\hat{g} \in \FH$, the constant term in the expansion of
$f(\hat{g}, u\- \cdot \bfu)$ is $f(\hat{g}, \bfu)$, and it is
invertible in $\mathscr{O}(\Omega(m))$, so $f(\hat{g}, u\- \cdot
\bfu)^{-1} $ belongs to $\mathscr{O}(\Omega(m))[[u]]$. Therefore, in
view of the expansion form (\ref{eq:psibfuexpansion}), the
coordinates of $\psi( u\- \cdot \bfu)$ expand into a formal power
series in $u_{\alpha}$ whose coefficients are series in
$\mathscr{O}(\Omega(m))$, but it follows from Lemma
\ref{lem:norms-m-m'} that the coefficients are indeed convergent series in
$\mathscr{O}(\Omega(m))$ for $u \in \RB $. So each coordinate of
$\psi(u\- \cdot \bfu)$ belongs to $\mathscr{O}(\Omega(m))[[u]]$.

For $l \in \RU^-_\RL $ or $\RT$,
$$\pi_\sigma (l) \psi(\bfu) = \sigma (l)\- \psi(l\- \cdot \bfu \cdot l).$$
The arguments are similar.
\end{proof}

\begin{cor}\label{cor:f(u+bfu)-1} Let $\RU^+_\fo = \RU^+ \cap
\RG_\fo$, then the power series expansion of $$f(j(u^+ , \bfu))^{-1}
= f(u^+ \cdot \bfu)\-, \hskip 10pt u^+ \in \RU^+_\fo,$$ on
$\RU^+_\fo$ converges in $\mathscr{O}(\Omega(m))$.
\end{cor}

\begin{proof}
Since $f$ is an $F$-rational character on $\RP^+$ (Lemma \ref{lem:f(g)} (1)), if we
put $\sigma = f$ and $\psi \equiv 1$, then $( \pi_{f} (g\-) 1) (\bfu) = f(j(g , \bfu))\-$.
Therefore our assertion follows from Lemma \ref{lem:pi_sigma analytic}.
\end{proof}

Now consider the dual representation $\pi_\sigma^*$ of $\RG$ on
$\mathscr{O}_\sigma(\Omega)^*_b \cong
\underset{m}{\varinjlim}\mathscr{O}(\Omega(m))^*_b$. The transition
homomorphisms $\mathscr{O}_\sigma(\Omega(m))^*_b \ra
\mathscr{O}_\sigma(\Omega)^*_b$ are injective (see Proposition
\ref{prop:Stein} (2)). Lemma \ref{lem:gOmegam} implies that, for any
$g\in \RG_\fo$, $\pi_\sigma^*(g)$ maps $\mathscr{O}(\Omega(m))^*_b$
into $\mathscr{O}(\Omega(m'))^*_b$ via
$$\left\langle\psi,~ \pi_\sigma^*(g) \mu \right\rangle = \left \langle (\pi_\sigma(g^{-1})\psi)|_{\Omega(m)}, ~ \mu \right \rangle,
\hskip 10pt \mu \in \mathscr{O}_\sigma(\Omega(m))^*, \psi \in
\mathscr{O}_\sigma(\Omega(m')).$$ We deduce from Lemma
\ref{lem:pi_sigma analytic} that, for any $\mu \in
\mathscr{O}_\sigma(\Omega(m))^*$, the orbit map
\begin{eqnarray*}  \RB ^{-1}  & \rightarrow & \mathscr{O}_\sigma(\Omega(m'))^*_b  \\
g \quad & \mapsto & \pi_{\sigma}^*(g) \mu
\end{eqnarray*}
is an $\mathscr{O}_\sigma(\Omega(m'))^*_b$-valued analytic function.
Therefore we obtain the following corollary.
\begin{cor}\label{cor:dual locally analytic}
$(\mathscr{O}_\sigma(\Omega)^*_b, \pi_\sigma^*)$ is locally
analytic.
\end{cor}

\section{Duality}\label{sec:Duality}

In the following, we assume that $K$ is spherically complete. Let
$(V, \sigma)$ be a $d$-dimensional $K$-rational representation of
$\RL$. We choose a basis $v_1, \cdots, v_d$ of $V$ and denote by
$v_1^*, \cdots, v_d^*$ the corresponding dual basis of the dual
space $V^*$. $(V^*, \sigma^*)$ denotes the dual representation of
$(V, \sigma)$.

\subsection{The duality operator $I_{\sigma}$} \label{sec:Isigma}

For $\bfu \in \Omega$ and $v^*\in V^*$, let $\varphi_{\bfu, v^*}$ be
the $V^*$-valued locally analytic function on $\FH$:
\begin{equation*}\label{eq:varphibfuv*}\varphi_{\bfu, v^*}(\hat{g}):=\sigma^*(j(\hat{g}, \bfu))v^*.
\end{equation*}

In view of (\ref{eq:j(lg,bfu)}), $\varphi_{\bfu, v^*}$ belongs to
$C^\an_{\sigma^*}(\FH, V^*)$. Let $B^0_{\sigma^*}(\FH, V^*)$ be the
subspace of $C^\an_{\sigma^*}(\FH, V^*)$ spanned by $\varphi_{\bfu,
v^*}$, $B_{\sigma^*}(\FH, V^*)$ the closure of $B^0_{\sigma^*}(\FH,
V^*)$. From (\ref{eq:automorphyfactor}), we see that $B^0
_{\sigma^*}(\FH, V^*)$ and therefore $B_{\sigma^*}(\FH, V^*)$ are
$\RG$-invariant.

For any continuous linear functional $\xi\in B _{\sigma^*}(\FH,
V^*)^* $, we define a $V$-valued function on $\Omega$:
\begin{equation*} \label{eq:Isigma}
I_\sigma(\xi) (\bfu) := \sum_{k=1}^d \langle \varphi_{\bfu, v^*_k},
\xi \rangle v_k, \hskip 10pt \bfu\in \Omega.
\end{equation*}
 $I_\sigma(\xi)$ is independent of the choice of the basis  $\{v_k
\}_{k=1}^d$. Evidently, $I_\sigma$ is injective.

\begin{lem} \label{lem:Iequivariant}
$I_\sigma$ is $\RG$-equivariant, that is,
$$I_\sigma(T_{\sigma^*}^*(g) \xi) = \pi_\sigma (g) I_\sigma (\xi),$$
for any $g\in \RG$.
\end{lem}
\begin{proof}
\begin{align*} 
 I_\sigma(T_{\sigma^*}^*(g) \xi)(\bfu) 
& =  \sum_{k=1}^d  \langle \varphi_{\bfu, v_k^*}, T_{\sigma^*}^*(g) \xi \rangle v_k 
 =  \sum_{k=1}^d  \langle T_{\sigma^*}(g^{-1})\varphi_{\bfu, v_k^*}, \xi \rangle v_k \\
&= \sum_{k=1}^d  \langle \sigma^*( j(? \cdot g\-, \bfu))v^*_k, \xi \rangle v_k  \\
& = \sigma (j(g^{-1}, \bfu))^{-1} \left(\sum_{k=1}^d  \langle
\sigma^*( j(?, g\- * \bfu) ) v_{k; g} ^*, \xi \rangle v_{k; g}\right) \hskip 10pt (\text{see } (\ref{eq:automorphyfactor}))\\
&= (\pi_\sigma (g) I_\sigma (\xi))(\bfu), 
\end{align*}
where $v_{k; g} = \sigma (j(g^{-1}, \bfu))v_k$, and similarly $ v_{k; g} ^* = \sigma^* (j(g^{-1}, \bfu))v^*_k$.
$\{v_{k; g} \}_{k=1}^d$ and $\{v_{k; g}^* \}_{k=1}^d$ are dual to each other.
\end{proof}

\begin{prop}\label{prop:Isigma}$\ $
	
$(1)$ For any continuous linear functional $\xi\in B _{\sigma^*}(\FH,
V^*)^* $, $I_\sigma(\xi)$ is a $V$-valued rigid analytic function on
$\Omega$. 

$(2)$ $ I_\sigma $ is a continuous homomorphism of
$\RG$-representations from $(B _{\sigma^*}(\FH, V^*)^*_b,
T^*_{\sigma^*})$ to $(\mathscr{O}_\sigma(\Omega), \pi _{\sigma})$.
\end{prop}
\begin{proof}

\emph{Step 1.} We denote by $i$ the inclusion: $B _{\sigma^*}(\FH,
V^*) \hookrightarrow C^\an_{\sigma^*}(\FH , V^*)$, $i^*$ its adjoint
operator. Because of our assumption that $K$ is spherically
complete, the Hahn-Banach Theorem (\cite {Schneider} Corollary 9.4)
implies that $i^*$ is surjective. Since $C^\an_{\sigma^*}(\FH ,
V^*)^*_b$ and $B _{\sigma^*}(\FH, V^*)^*_b$ are both Fr\'echet
spaces (Corollary \ref{cor:compacttype}), the
open mapping theorem (\cite{Schneider} Proposition 8.6) implies that $i^*$ is open.
Therefore the
continuity of $I_\sigma\circ i^*$ implies that of $I_\sigma$.
Consequently, (1) and (2) are equivalent to:\\
($1'$) $I_\sigma\circ i^*(\xi)\in \mathscr{O}_\sigma(\Omega)$ for
any $\xi\in C^\an_{\sigma^*}(\FH, V^*)^* $;
\\
($2'$) $I_\sigma\circ i^*: (C^\an_{\sigma^*}(\FH , V^*)^*_b,
T^*_{\sigma^*}) \ra (\mathscr{O}_\sigma(\Omega), \pi _{\sigma})$ is
a continuous homomorphism of $\RG$-representations.

Since $\RG$-equivariance is proved in Lemma \ref{lem:Iequivariant},
for ($2'$) it remains to show the continuity of $I_\sigma\circ i^*$.

For convenience, we still denote $I_\sigma\circ i^*$ by $I_\sigma$.

\emph{Step 2. } Let $\{\overline{\FU}_{\kappa }\}_{\kappa} $ be a
finite disjoint open covering
of ${\overline{\FH}} $ satisfying: \\
1. $\RU^+_\fo \in\{\overline{\FU}_{\kappa }\}_{\kappa} $
(note that the open subscheme $\RP^-\backslash \RC $ of $\overline{\FH}$ is identified with $\RU^+$); \\
2. each $\overline{\FU}_{\kappa }$ is (right) translated into
$\RU^+_\fo$ by
some $g_\kappa\in \RG$.\\
 Let $\FU_{\kappa }$ be the preimage of
$\overline{\FU}_{\kappa }$ under $
 \mathrm{pr}^{\FH}_{\overline{\FH}} $.


For $\xi\in C^\an_{\sigma^*}(\FH, V^*)^* $, we write $I_\sigma(\xi)$
in integral:
\begin{align*} 
I_{\sigma}(\xi)(\bfu) =  \sum_{k=1}^d \int_{\FH} \varphi_{\bfu; v^*_k} d \xi\cdot v_k 
& =  \sum_{k=1}^d \sum_{\kappa } \int_{\FU_{\kappa }}
\varphi_{\bfu; v^*_k} d \xi\cdot v_k\\
& = \sum_{\kappa } \pi_{\sigma}(g_{\kappa }) \Big( \sum_{k=1}^d
\int_{\FU_{\kappa }\cdot g_{\kappa }} \varphi_{\bfu; v_{k; g_\kappa
} ^*} d(T^*_{\sigma^*}(g_{\kappa }^{-1})\xi)\cdot v_{k; g_\kappa
}\Big),
\end{align*}
where $v_{k; g_\kappa} = \sigma (j(g_\kappa^{-1}, \bfu))v_k$ is
defined in the proof of Lemma \ref{lem:Iequivariant}. Therefore it
suffices to consider
\begin{equation}\label{eq:integerU}
\sum_{k=1}^d \int_{\FU} \varphi_{\bfu; v^*_k} d \xi'\cdot v_k,
\end{equation}
where  $\FU$ ranges on $\{\FU_{\kappa }\cdot g_{\kappa } \}_\kappa $ and
$\xi'$ is the image of $\xi$ under $C^\an_{\sigma^*}(\FH, V^*)^*_b
\ra C^\an_{\sigma^*}(\FU, V^*)^*_b$.

For the open subset $\overline{\FU} =
\mathrm{pr}^\FH_{\overline{\FH}} (\FU)$ of $\RU^+_\fo$, we have the
isomorphism induced from a locally analytic section $\iota$ of
$\mathrm{pr}^\FU_{\overline{\FU}}$ (see Lemma
\ref{lem:principalseries} (3)):
\begin{equation}\label{eq:isoCO} C^\an_{\sigma^*}(\FU , V^*)^*_b\simeq C^\an(\overline{\FU}, V^*)^*_b.
\end{equation}
Then (\ref{eq:integerU}) is equal to $$\bar{I}_{\sigma,
\overline{\FU}}(\overline{\xi})(\bfu) := \sum_{k=1}^d
\int_{\overline{\FU}} (\sigma^*(j(u^+, \bfu))v^*_k)\ d
\overline{\xi}(u^+)\cdot v_k,$$ where $\overline{\xi}$ is the image
of $\xi'$ in $C^\an(\overline{\FU}, V^*)^*_b$ via the isomorphism
(\ref{eq:isoCO}).

Therefore it suffices to prove that $\bar{I}_{\sigma,
\overline{\FU}}(\overline{\xi})$ is rigid analytic on $\Omega(m)$,
and that the map
\begin{eqnarray*}C^\an(\overline{\FU}, V^*)^*_b &\ra& \mathscr{O}_\sigma(\Omega(m))\\
\overline{\xi} &\mapsto& \bar{I}_{\sigma,
\overline{\FU}}(\overline{\xi})|_{\Omega(m)}
\end{eqnarray*}
is continuous for all $m$.

\emph{Step 3.} Since $\sigma^*$ is $K$-rational, using the same
arguments in the proof of Lemma \ref{lem:pi(psi)inO} and applying Corollary
\ref{cor:f(u+bfu)-1}, we obtain an expansion $$ \sigma^*(j(u^+,
\bfu)) v_k^*=\sum_{\ell=1}^d \Big( \sum_{{\underline{r}}\in
\BN_0^{R_I^+}}a_{{\underline{r}}, k \ell} (\bfu)\cdot
(\underline{u}^+)^{\underline{r}} \Big) v^*_\ell,$$
 with $a_{{\underline{r}}, k \ell}\in \mathscr{O}(\Omega(m))$ such
 that
\begin{equation}\label{eq:a_r,kl 1}\underset{|{\underline{r}}|\ra
\infty}{\lim}||a_{{\underline{r}}, k
\ell}||_{\mathscr{O}(\Omega(m))} \cdot
||(\underline{u}^+)^{\underline{r}}||_{C^\an(\RU^+_\fo)} = 0,
\end{equation}
and moreover, there is a constant $c'>0$, depending only on $m$,
$\sigma$ and $\{v_k\}_{k=1}^d$, such that
\begin{equation}\label{eq:a_r,kl 2} ||a_{{\underline{r}}, k \ell}||_{\mathscr{O}(\Omega(m))}
\cdot ||(\underline{u}^+)^{\underline{r}}||_{C^\an(\RU^+_\fo)} \leqslant 
c'.
\end{equation}

Then
\begin{equation} \label{eq:expansionbarI} \bar{I}_{\sigma, \overline{\FU}}(\overline{\xi})(\bfu) =\sum_{k=1}^d \Big(\sum_{\ell=1}^d
\sum_{{\underline{r}}} \int_{\overline{\FU}}
 (\underline{u}^+)^{\underline{r}}\cdot v^*_\ell d\overline{\xi}(u^+) \cdot a_{{\underline{r}}, k\ell } (\bfu) \Big) v_k.
\end{equation}

We have
\begin{equation}\label{eq:intbfunvl} \left|\int_{\overline{\FU}}
(\underline{u}^+)^{\underline{r}}\cdot v^*_\ell
d\overline{\xi}(u^+)\right|\leqslant 
||(\underline{u}^+)^{\underline{r}}||_{C^\an(\RU^+_\fo)} \cdot
||v^*_\ell||_{V^*} \cdot ||\overline{\xi}||_{C^\an(\overline{\FU},
V^*)^*_b}.\end{equation}

(\ref{eq:a_r,kl 1}) and (\ref{eq:intbfunvl}) imply that the
expansion (\ref{eq:expansionbarI}) of $\bar{I}_{\sigma,
\overline{\FU}}(\overline{\xi})$ converges in
$\mathscr{O}_\sigma(\Omega(m))$.

(\ref{eq:a_r,kl 2}) and (\ref{eq:intbfunvl}) imply
$$\left\| \bar{I}_{\sigma, \overline{\FU}}(\overline{\xi})\right\|_{\mathscr{O}_\sigma(\Omega(m))}\leqslant 
\max_{1\leqslant  k, \ell\leqslant  d} c' ||v_\ell^*||_{V^*} ||v_k||_V\cdot
||\overline{\xi}||_{C^\an(\overline{\FU}, V^*)^*_b},
$$ and therefore the continuity follows.
\end{proof}

\subsection{The duality operator $J_\sigma$}

Let $\mathscr{N}_\sigma (\Omega)$ denote the image of $I_\sigma$.

We consider $J_\sigma$, the adjoint operator of $I_\sigma$, which is
an injective continuous linear operator from 
$\mathscr{N}_\sigma(\Omega)^*_b$ to $(B_{\sigma^*}(\FH,
V^*)^*_b)_b^{*}\cong B_{\sigma^*}(\FH, V^*)$ ($B_{\sigma^*}(\FH,
V^*)$ is reflexive according to Corollary \ref{cor:compacttype}).

For any $\mu\in \mathscr{N}_\sigma(\Omega)^*$ and $\xi\in
B_{\sigma^*}(\FH, V^*) ^*$, we have 
\begin{equation} \label{eq:dual} \langle J_\sigma(\mu), \xi\rangle=\langle I_\sigma(\xi), \mu\rangle.
\end{equation}

For $\hat{g}\in \FH$ and $v\in V$, we define the Dirac distribution
$\xi_{\hat{g}, v} \in B_{\sigma^*}(\FH, V^*) ^*$ as follows:
$$\langle \varphi, \xi_{\hat{g}, v}\rangle = \langle v,
\varphi(\hat{g})\rangle_V, \hskip 10pt \varphi \in B_{\sigma^*}(\FH,
V^*), $$ and a $V$-valued rigid analytic function $\psi_{\hat{g},
v}$ on $\Omega$: $$\psi_{\hat{g}, v} (\bfu) := \sigma(j(\hat{g},
\bfu))^{-1}v.$$

\begin{lem}\label{lem:Jsigma-xi} $$I_\sigma(\xi_{\hat{g},
v})=\psi_{\hat{g}, v} . $$
\end{lem}
\begin{proof} This is straightforward from definitions.
	\begin{align*}
	I_\sigma(\xi_{\hat{g},		v}) (\bfu) 
	  = \sum_{k=1}^d \langle \varphi_{\bfu, v^*_k},
	\xi_{\hat{g}, v} \rangle v_k 
	& = \sum_{k=1}^d \langle v,  \varphi_{\bfu, v^*_k} (\hat{g}) \rangle_V v_k 
	  = \sum_{k=1}^d \langle v, \sigma^*(j(\hat{g}, \bfu))v^*_k  \rangle_V v_k \\
	& = \sum_{k=1}^d  \langle \sigma (j(\hat{g}, \bfu)) v, v^*_k  \rangle_V v_k = \sigma (j(\hat{g}, \bfu)) v = \psi_{\hat{g}, v} (\bfu)
	\end{align*}
\end{proof}

Then we obtain a formula for $J_{\sigma}$.

\begin{prop} For any continuous linear
functional $\mu\in \mathscr{N}_\sigma(\Omega)^*$, we have
\begin{equation}\label{eq:Jsigma}
J_{\sigma}(\mu)(\hat{g})=\sum_{k=1}^d \langle\psi_{\hat{g}, v_k},
\mu\rangle v_k^*.
\end{equation}
\end{prop}
\begin{proof} This is straightforward from Lemma \ref{lem:Jsigma-xi} and (\ref{eq:dual}). Indeed,
{\allowdisplaybreaks
		\begin{align*} \sum_{k=1}^r  \langle\psi_{\hat{g}, v_k},
			\mu\rangle v_k^* = \sum_{k=1}^r \langle I_\sigma(\xi_{\hat{g},
				v_k}), \mu
			\rangle v_k^*   
			   & = \sum_{k=1}^r\langle J_\sigma(\mu), \xi_{\hat{g}, v_k}\rangle
			v_k^* \\ 
			& =  \sum_{k=1}^r \langle v_k, J_\sigma(\mu)(\hat{g}) \rangle_V
			\cdot v_k^* = 
			J_\sigma(\mu)(\hat{g}).
		\end{align*} }
\end{proof}

\subsection{The image of $I_\sigma$}

Let $\mathscr{N}_\sigma^0(\Omega)$ denote the subspace of
$\mathscr{O}_\sigma(\Omega)$ spanned by $\psi_{\hat{g}, v}$ for all
$\hat{g}\in \FH $ and $v\in V$. Then it follows from
(\ref{eq:automorphyfactor}) that $\mathscr{N}^0_\sigma(\Omega)$ is
$\RG$-invariant, and Lemma \ref{lem:Jsigma-xi} implies
$\mathscr{N}_\sigma^0(\Omega)\subset \mathscr{N}_\sigma(\Omega)$.

From (\ref{eq:Jsigma}), we see that $J_\sigma$ factors through
$\mathscr{N}_\sigma^0(\Omega)^*$, and (\ref{eq:Jsigma}) defines an
injective map from $\mathscr{N}^0_\sigma(\Omega)^*_b $ into $
B_{\sigma^*}(\FH, V^*)$. Since $J_\sigma$ is injective and the
natural map $\mathscr{N}_\sigma(\Omega)^*_b\ra
\mathscr{N}^0_\sigma(\Omega)^*_b$ is surjective (the Hahn-Banach
Theorem), $\mathscr{N}^0_\sigma(\Omega)^*_b =
\mathscr{N}_\sigma(\Omega)^*_b$. Therefore the Hahn-Banach Theorem
implies the following lemma.

\begin{lem}\label{lem:CN0denseinCN}$\mathscr{N}^0_\sigma(\Omega)$ is dense in $\mathscr{N}_\sigma(\Omega)$.
\end{lem}

\begin{thm} \label{thm:Iisomorphism}~ \\
 $(1)$ $I_\sigma $ is an isomorphism from $B _{\sigma^*}(\FH, V^*)^*_b$ to
$\mathscr{N}_\sigma(\Omega)$. \\
 $(2)$ $\mathscr{N}_\sigma(\Omega)$ is the
closure of $\mathscr{N}^0_\sigma(\Omega)$ in
$\mathscr{O}_\sigma(\Omega)$.
\end{thm}
\begin{proof}

 Let $\iota$ be a locally analytic section of
$\mathrm{pr}^\FH_{{\overline{\FH}}}$, and denote $\FK =
\iota(\overline{\FH})$.

\emph{1.} Let $\mathscr{N}^0_\sigma(\Omega(m))$ be the image of
$\mathscr{N}^0_\sigma (\Omega )$ in $\mathscr{O}_\sigma(\Omega(m))$.

Since $\psi_{\hat{g}, v_k} = \pi_\sigma(g\-) v_k$, we see that the
map
\begin{eqnarray*} \FH & \ra & \mathscr{O}_\sigma(\Omega(m)) \\
\hat{g} &\mapsto & \psi_{\hat{g}, v_k},
\end{eqnarray*}
 is locally analytic (Lemma \ref{lem:pi_sigma analytic}). Since $\FK
$ is compact, $\rho_m = \underset {1\leqslant  k\leqslant  d} {\min}\underset
{\hat{g}\in \FK} {\min} ||\psi_{\hat{g},
v_k}||_{\mathscr{O}_\sigma(\Omega(m))}$ is positive. Let
$\mathscr{L}$ denote the lattice $\ds \sum_{k=1}^d
\underset{\hat{g}\in \FK}{\sum} \fo_K\cdot \psi_{\hat{g}, v_k} $ in
$\mathscr{N}^0_\sigma(\Omega)$. Then, for each $m$, the image of
$\mathscr{L}$ in $\mathscr{N}^0_\sigma(\Omega(m))$ contains the ball
of radius $\rho_m$ centered at zero, and therefore the interior of
$\mathscr{L}$ is a nontrivial open lattice.

\emph{2. } According to Lemma \ref{lem:principalseries}, $\iota $ induces an isomorphism
$\iota^\circ$ between $C^\an_{\sigma^*}(\FH , V^*)$ and
$C^\an(\overline{\FH}, V^*)$,
and hence an isomorphism between $B_{\sigma^*}(\FH, V^*)$ and its image,  denoted by $B(\overline{\FH}, V^*).$

Let $\CI$ be any (finite) disjoint open chart covering
$\{\overline{\FU}_\kappa\}_{\kappa}$ of $\overline{\FH}$. We
recall that $C^\an(\overline{\FH}, V^*)$ is defined to be the inductive
limit, indexed with all the $\CI$, of the $K$-Banach algebras
$E_{\CI}(\overline{\FH}, V^*) = \prod_{\kappa} \CO
(\overline{\FU}_\kappa, V^*)$, where $\CO (\overline{\FU}_\kappa,
V^*)$ denotes the space of $K$-analytic functions on
$\overline{\FU}_\kappa$ (cf. \cite{Feaux} 2.1.10 and
\cite{Schneider-Teitelbaum2002} \S 2). The inductive limit structure
is naturally induced onto $B(\overline{\FH}, V^*)$, say
$B(\overline{\FH}, V^*) = {\varinjlim}_{\CI} E_{\CI}(\overline{\FH},
V^*).$ Moreover, the strong dual space $B(\overline{\FH}, V^*)^*_b$ is the
projective limit of $E_{\CI}(\overline{\FH}, V^*)^*_b.$

\emph{3. }Consider
\begin{eqnarray*}  \big({\iota ^{\circ}}\- \big)^*  \circ I_\sigma\- |_{\mathscr{N}^0_\sigma(\Omega )}:
\mathscr{N}^0_\sigma(\Omega )& \ra & B (\overline{\FH}, V^*)^*_b\\
\psi_{\hat{g}, v} &\mapsto& \big ({\iota ^{\circ}}\- \big)^*(\xi_{\hat{g},
v}).
\end{eqnarray*}

Let $\hat{g}\in \FK$.
\begin{align*} 
\left\| \big({\iota ^{\circ}}\- \big)^* (\xi_{\hat{g}, v})\right\|_{E_{\CI} (\overline{\FH}, V^*)^*_b}& =
\max _{\overline{\varphi}\in E_{\CI} (\overline{\FH}, V^*)}
\frac{\langle \overline{\varphi} ,\ \big({\iota ^{\circ}}\-\big)^*
(\xi_{\hat{g}, v})
 \rangle }{||\overline{\varphi}||_{E_{\CI} (\overline{\FH}, V^*)}} \\
& =  \max _{\varphi \in\ {\iota ^{\circ}}\- (E_{\CI}
(\overline{\FH}, V^*))}\frac{\langle \varphi ,\ \xi_{\hat{g}, v}
 \rangle}{||\iota ^{\circ}(\varphi)||_{E_{\CI} (\overline{\FH}, V^*)}}\\
& =  \underset{\varphi\in\ {\iota ^{\circ}}\- (E_{\CI}
(\overline{\FH}, V^*))} {\max}\frac{\langle v,
\varphi(\hat{g})\rangle_V}{\underset{\hat{g'}\in \FK}
{\max}||\varphi(\hat{g'})||_{V^*}}\\
&\leqslant   ||v||_V.
\end{align*}

Therefore the image of $\mathscr{L}$ under $\big({\iota ^{\circ}}\- \big)^*
\circ I_\sigma\- |_{\mathscr{N}^0_\sigma(\Omega )}$ in $B
(\overline{\FH}, V^*)^*_b$ is bounded, since its image in $E_{\CI}
(\overline{\FH}, V^*)^*_b$ are all norm-bounded by $\underset{1\leqslant 
k\leqslant  d}{\max}||v_k||_V$. Because $\mathscr{N}^0_\sigma(\Omega)$ is
metrizable, it is bornological (\cite{Schneider} Proposition 6.14),
and therefore $I_\sigma^{-1}|_{\mathscr{N}^0_\sigma(\Omega )}$ is
continuous (\cite{Schneider} Proposition 6.13). Therefore $I_\sigma$
induces an isomorphism between $I_\sigma\-
(\mathscr{N}^0_\sigma(\Omega))$ and $\mathscr{N}^0_\sigma(\Omega)$, and consequently $I_\sigma$ induces an isomorphism between their completions,
which, in view of Lemma \ref{lem:CN0denseinCN}, must be $B
_{\sigma^*}(\FH, V^*)^*_b$ and $\mathscr{N}_\sigma(\Omega)$ respectively.
\end{proof}




\delete{ By (\ref{eq:Jsigma}), {\allowdisplaybreaks
\begin{eqnarray*}
\big(I_\sigma(\m\CU_{\bfu,v^*})\big)(\hat{g}) &=&
\sum_{k=1}^r\langle \psi_{\hat{g}, v_k}, \m\CU_{\bfu,v^*}\rangle v_k^*  \\
&=& \sum_{k=1}^r \langle\sigma(j(\hat{g}, \bfu))^{-1}v_k, v^*\rangle_V \cdot v_k^* \\
&=& \sum_{k=1}^r \langle v_k, \sigma^*(j(\hat{g}, \bfu)) v^*\rangle_V \cdot v_k^*  \\
&=& \sigma^*(j(\hat{g}, \bfu)) v^* \\
&=& \varphi_{\bfu,v^*}(\hat{g}).
\end{eqnarray*}
} }


\begin{cor}
$J_\sigma$ is an isomorphism of $\RG$-representations from
$(\mathscr{N}_\sigma(\Omega)^*_b, \pi_{\sigma}^*)$ to
$(B_{\sigma^*}(\FH, V^*), T_{\sigma^*})$.
\end{cor}

\section{Concluding remarks}

In \cite{QiYang} \S 3 we briefly reviewed Morita's theory of $\SL(2,
F)$ and discussed the relation between $I_\sigma$ and Morita's
duality and Casselman's operator for $$\sigma_s\begin{pmatrix}z\- & 0
\\ 0 & z \end{pmatrix} = z^s$$ with $s$ a positive integer
 (for $s$ non-positive, $I_\sigma$ is an isomorphism between two
$(-s + 1)$-dimensional $G$-representations, which is of less
interest).

To illustrate this connection, we consider the special case $s = 2$. $\mathscr{O}_{ \sigma_2} (\Omega) $ is canonically isomorphic to
the space $\Omega^1(\Omega)$ of holomorphic $1$-forms on the upper half plane $\Omega \cong K - F = \BP^1(K) - \BP^1(F)$ 
via $\psi(\bfu) \mapsto \psi(\bfu) d \bfu.$ It may be shown that $\mathscr {N}_{\sigma_2} (\Omega)$ corresponds
to the subspace of $\Omega^1(\Omega)$ with zero residue at each point of $ \BP^1(F)$ (see \cite{MoritaI} and \cite {QiYang}). On the other hand,
$ C^{\an}_{\sigma_0}(\FH) \cong C^{\an} (\overline {\FH})$ with $\overline {\FH} \cong \BP^1(F),$ and we denote $D_0 = C^{\an} (\BP^1(F))$. $D_0$ has two closed $\RG$-invariant subspaces, the spaces $P_0$ and $P^{\loc}_0$ consisting of constants and locally constant functions on $\BP^1(F)$ respectively. The classical Morita's duality is established via residues. More precisely, for each $\psi \in \mathscr{O}_{ \sigma_2} (\Omega)$, define a linear functional $  M  _2(\psi)$ of $D_0$ by
$$ \left \langle \varphi ,   M  _2(\psi) \right \rangle =  \text{ the sum of residues
of the } 1 \text{-form } \varphi(u) \psi(u)\ du \text{ on } \BP^1(F).$$
Morita's duality $  M  _2$ induces $\RG$-isomorphisms $\mathscr{O}_{ \sigma_2} (\Omega) \cong (D_0/P_0)^*_b $ and $\mathscr {N}_{\sigma_2} (\Omega) \cong (D_0/P^{\loc}_0)^*_b$. Moreover, Casselman's intertwining operator $$S_0: \varphi \mapsto  d \varphi $$ induces a $\RG$-isomorphism between $D_0/P_0^{\loc}$ and the space $D_{-2}$ of locally analytic $1$-forms on $\BP^1(F)$, a space that is isomorphic to $ C^{\an}_{\sigma_{-2}}(\FH)$ (see \cite{MoritaII} and \cite {QiYang}). The connection between our duality operator $I_{\sigma_{2}}$ and Morita's duality $  M  _2$ was found in 
\cite {QiYang} Theorem 3.6 as the following commutative diagram
\begin{center}
\begin{tikzpicture}
\matrix(m)[matrix of math nodes,
row sep=3em, column sep=0em,
text height=1.5ex, text depth=0.25ex]
{\mathscr{N}_{\sigma_2}(\Omega) & &(D_{0}/P^\loc_{0})^*_b\\
B_{\sigma_{-2}}(\FH )^*_b &= \hskip 2.5pt  C^{\an}_{\sigma_{-2}}(\FH)^*_b &  \cong \hskip 2.5pt (D_{-2})^*_b \\};
\path[->,font=\scriptsize]
(m-1-1) edge node[auto] {$  M  _2$} (m-1-3)
(m-2-1) edge node[auto] {$I_{\sigma_2}$} (m-1-1)
(m-2-3) edge node[auto] {$S_{0}^*$} (m-1-3);
\end{tikzpicture}
\end{center}

A generalization of Morita's duality seems quite hard in view of its
analytic construction via residues. The first step towards this
would be finding other closed sub-representations of
$(C^\an_\sigma({\FH} , V), T_{\sigma})$ and $(
\mathscr{O}_\sigma(\Omega), \pi_{\sigma})$. This work is done
completely for $\SL(2, F)$ in Morita and Murase's \cite{MoritaI},
\cite{MoritaII} and \cite{MoritaIII}, where the complete
classifications of the sub-quotient spaces of holomorphic discrete
series and the principal series are conjectured and claimed (Morita
attempted to prove this, but his proof contained a serious gap).

A further question is on the irreducibility. For this, we conjecture that
$(\mathscr{N}_\sigma(\Omega), \pi_\sigma)$ and $(B_{\sigma^*}(\FH,
V^*) , T_{\sigma^*})$ are topologically irreducible
$\RG$-representations if $\sigma$ is irreducible. For $\SL(2, F)$,
this conjecture was claimed in \cite{MoritaIII} Theorem 1 (i) and a
proof for $F = \BQ_p$ was given by Schneider and Teitelbaum in
\cite{Schneider-Teitelbaum2002}.

\end{document}